%% file: paper.tex
\documentclass{amsart}

\usepackage{amsmath,amsthm,amssymb}
\usepackage{mathtools}
\usepackage{graphicx}
\usepackage{cite}
\usepackage{hyperref}
\usepackage[capitalise]{cleveref}
\Crefname{equation}{}{}

\newtheorem{theorem}{Theorem}
\newtheorem{lemma}[theorem]{Lemma}
\newtheorem{proposition}[theorem]{Proposition}
\newtheorem{corollary}[theorem]{Corollary}

\input{macros}

\author{Ansgar Freyer}
\address{Institut f\"ur Diskrete Mathematik und Geometrie,
Technische Universit\"at Wien,
Wiedner Hauptstra\ss e 8-10/1046,
1040 Wien, Austria}
\email{ansgar.freyer@tuwien.ac.at}

\author{Monika Ludwig}
\address{Institut f\"ur Diskrete Mathematik und Geometrie,
Technische Universit\"at Wien,
Wiedner Hauptstra\ss e 8-10/1046,
1040 Wien, Austria}
\email{monika.ludwig@tuwien.ac.at}

\author{Martin Rubey}
\address{Institut f\"ur Diskrete Mathematik und Geometrie,
Technische Universit\"at Wien,
Wiedner Hauptstra\ss e 8-10/1046,
1040 Wien, Austria}
\email{martin.rubey@tuwien.ac.at}

\title{Unimodular Valuations beyond Ehrhart}

\begin{document}

\begin{abstract}
A complete classification of unimodular valuations on the set of lattice polygons with values in the spaces of polynomials and formal power series, respectively, is established. The valuations are classified in terms of their behaviour with respect to dilation using extensions to unbounded polyhedra and basic invariant theory.

\bigskip

\end{abstract}

\maketitle
\section{Introduction}

For $n\ge 2$, let $\lpn$ denote the set of lattice polytopes in $\RR^n$, that is,  convex polytopes with vertices in the integer lattice $\ZZ^n$. A function $\oZ$ on $\lpn$ with values in a vector space $V$ is a \Dfn{valuation}, if
\[
 \oZ(P)+\oZ(Q) =  \oZ(P\cup Q)+\oZ(P\cap Q)
\]
for all $P,Q\in \lpn$ such that both $P\cup Q$ and $P\cap Q$ are also in $\lpn$.  
The lattice point enumerator $\Lat^0\colon\lpn\to \RR$, which counts the number of lattice points in $P\in\lpn$, and $n$-dimensional volume are important valuations on $\lpn$.

A fundamental result on valuations on lattice polytopes is the Betke--Kneser Theorem  \cite{Betke:Kneser}. It provides a complete classification of real-valued valuations on lattice polytopes in $\RR^n$ that are invariant with respect to the \Dfn{affine unimodular group}, obtained by the general linear group over the integers $\glnz$ and translations by vectors from $\ZZ^n$. It also provides a characterization of the so-called Ehrhart coefficients \cite{Ehrhart62}. We state the two-dimensional version.

\begin{theorem}[Betke--Kneser]
\label{thm:betke-kneser}
A function $\,\oZ\colon\PZ\to \RR$ is  a translation and $\glz$ invariant valuation if and only if  there are $\,
c_0, c_1, c_2\in\RR$ such that 
$$\oZ(P) =  c_0  \Lat^0_0(P)+c_1 \Lat^0_1(P) +c_2 \Lat^0_2(P)$$
for every $P\in\PZ$. 
\end{theorem}

\noindent
In the planar case, the Ehrhart coefficients have a simple description: $\Lat^0_0(P)\coloneq 1$  for a non-empty lattice polygon $P$, while $\Lat^0_1(P)$ is half the number of lattice points in the boundary of $P$, and $\Lat^0_2(P)$ is the area of $P$.
We refer to \cite{Barvinok2008,BeckRobins,Gruber,Jochemko19, McMullen} for general information on valuations on lattice polytopes and to \cite{Bajo,  BergJochemkoSilverstein, BDHHL,  JochemkoSanyal17, JochemkoSanyal18, ludwigsilverstein, McMullen} for some recent results.
 
We generalize the above theorem and establish a complete classification of poly\-nomial-valued, unimodular valuations on $\PZ$. 
Here, we call a valuation \Dfn{uni\-modular} if it is translatively polynomial and $\glz$ equivariant, and we recall both definitions below. 

\goodbreak
Let $\ring$ denote the vector space of bivariate polynomials and consider the (left) action of $\GL_2(\ZZ)$ on $\ring$ defined by $\phi f\coloneq  f \circ \phi^\star$ for $f\in\ring$ and
$\phi\in\GL_2(\ZZ)$, where $\phi^\star$ is the transpose of $\phi$.
A function $\oZ\colon\PZ\to\ring$ is \Dfn{$\glz$ equivariant} if $\oZ(\phi P)= \phi \oZ(P)$ for every $\phi\in\glz$ and $P\in \PZ$.

\goodbreak
Since $\ring$ is the direct sum of the spaces of $r$-homogeneous polynomials $\ring_r$ for $r\in\ZZ_{\ge 0}$, it suffices to classify valuations with values in $\ring_r$. 
This corresponds to a classification of valuations from $\PZ$ to the space of symmetric \mbox{$r$-tensors}. In this setting, the question was first considered in \cite{ludwigsilverstein} (see also \cite[Chapter 21]{McMullen}), where a complete classification was established for $r\le 8$. We remark that $\slnz$ equivariance was considered in \cite{ludwigsilverstein}, but there is a mistake in a calculation in the planar case. However, the results are correct for $\glz$ equivariance in the planar case and for $\slnz$ equivariance for $n\ge3$. 

All tensor valuations that appear for $r\le 8$ are tensor-valued Ehrhart coefficients (introduced in \cite{ BoeroeczkyLudwig_survey, BoeroeczkyLudwig} for $r\ge 1$) obtained through a homogeneous decomposition of the $r$th \Dfn{discrete moment tensor}. This tensor corresponds to the $r$th discrete moment polynomial,
\[
\Lat^r(P) \coloneq  \frac{1}{r!} \sum_{v\in P\cap\ZZ^2} v^r,
\]
where  $v^r(x,y) = (ax + by)^r$ for $v=(a, b)$, that is, we consider $(x,y)$ as element of the dual of $\RR^2$ and $v$ as an element of its bidual. In particular, $\Lat^0$ is the classical lattice point enumerator. Similar to Ehrhart's celebrated result, there is a homogeneous decomposition of $\Lat^r$ (see \cite{BoeroeczkyLudwig_survey}): 
for $r\ge 1$ and $1\leq i\leq r+2$, there exist $i$-homogeneous valuations $\Lat_i^r\colon \PZ\to \ring_r$ such that 
\begin{equation}\label{eq:EC}
\Lat^r = \sum_{i=1}^{r+2}\Lat_i^r.
\end{equation}
Here, a function $\oZ$ is called \Dfn{$i$-homogeneous} if $\oZ(mP)=m^i\oZ(P)$ for every $m\in\ZZ_{\ge0}$ and $P\in\PZ$.

\goodbreak
A valuation $\oZ\colon \PZ\to\ring$ is \Dfn{translatively polynomial} (of degree at most $r$)  
if there exist \Dfn{associated valuations} $\oZ^{r-j}\colon\PZ\to\ring$ for $0\leq j \leq r$ such that
\begin{equation}
\label{eq:translation}
\oZ(P+v) = \sum_{j=0}^r \oZ^{r-j}(P) \frac{v^j}{j!}
\end{equation}
for every $P\in\PZ$ and $v\in\ZZ^2$. The associated valuations are uniquely determined, and by setting $v=0$, we see that $\oZ^r = \oZ$. 
Khovanski\u{i} and Pukhlikov \cite{pukhlikovkhovanskii} showed that every translatively polynomial valuation $\oZ\colon\PZ\to\ring$ of degree at most $r$ can be written as
\begin{equation*}
    \label{eq:hom_dec}
\oZ = \sum_{i=0}^{r+2} \oZ_i^r,
\end{equation*}
where $\oZ_i^r$ are $i$-homogeneous valuations.

Let $\Val$ denote the vector space of translatively polynomial and $\glz$ equivariant valuations $\oZ\colon\PZ\to \ring$. Let $\Val^r$ be the subspace of valuations $\oZ:\PZ\to \ring_r$. We call $r$ the \Dfn{rank} of the valuation.
It was proved in \cite{ludwigsilverstein} that the space
\[
\Val_i^r\coloneq \{\oZ\in\Val^r\colon \oZ \text{ is $i$-homogeneous}\}
\]
is spanned by the Ehrhart coefficient $\Lat^r_i$ for $1\le r \le 8$ and $1\le i\le r+2$, and that it is zero-dimensional for $i=0$ and $r\ge1$, as well as for $i>r+2$. Moreover, it was proved that $\dim \Val_1^9\ge 2$.

\goodbreak
Our first main result determines the dimensions of the vector spaces $\Val_i^r$. 

\begin{theorem}
\label{thm:tensor_val}
    For $i,r\in\ZZ_{\ge0}$, 
   \[
   \dim \Val_i^r= \begin{cases}
        1   &\text{ for } 1\le i<r \text{ and } r-i\text{ odd, and for }r\le i\le r+2,\\
        \parts{r}&\text{ for } i=1<r \text{ and } r-1 \text{ even},\\
        \parts{\frac{r-i}{2}+1} &\text{ for } 1<i<r \text{ and } r-i\text{ even},\\
         0 &\text{ otherwise, }
   \end{cases}\]
where $\parts{r} \coloneq  |\{(k,\ell)\in\ZZ_{\geq 0}^2 \colon 2k+3\ell=r\}|$
is the number of integer partitions of $r$ into parts $2$ and $3$.
\end{theorem}

\noindent
The vector space $\Val_i^r$ is spanned by $\Lat_i^r$ if $\dim  \Val_i^r=1$, and $\Lat_i^r\ne0$  (that is, $\Lat_i^r$ is not identically zero) for every $1\le i \le r+2$. Since $\parts{r}=1$ for $r\in\{3,5,7\}$, \Cref{thm:tensor_val} gives a basis of $\Val_1^r$ for every even $r$ and for $r\le 7$. 

The following result provides bases of the vector spaces $\Val_1^r$ also for the remaining~$r$. 
Let $\sttriangle$ be the standard triangle with vertices at the origin,  $e_1\coloneq (1,0)$, and $e_2\coloneq (0,1)$. Recall that every $P\in\PZ$ has a (generally non-unique) unimodular triangulation $\cT$; that is, $P$ can be written as finite union of non-overlapping triangles $\phi \sttriangle +v\in\cT$ with $\phi\in\glz$ and $v\in\ZZ^2$.

\begin{theorem}\label{thm:basis}
For $r>1$ odd, a basis of $\Val_1^r$ is given by $\Lat_1^{2k,3\ell}\colon\PZ\to\ring_r$ with $k,\ell\in \ZZ_{\ge0}$ and $2k+3\ell=r$, defined for $P\in\PZ$ and $\,\cT$  any unimodular triangulation of $P$, by
\begin{equation*}\label{eq:p23sum}
\Lat_1^{2k,3\ell}(P)\coloneq  \sum\nolimits_{\phi\sttriangle +v \in \cT} \big(\Lat_1^2(\sttriangle)^k \Lat_1^3(\sttriangle)^\ell\big)\circ \phi^\star
\end{equation*}
Here, $\Lat_1^{2k,3\ell}(P)$ does not depend on the triangulation $\cT$ of $P$.
\end{theorem}
\noindent
Note that it follows that the Ehrhart coefficient $\Lat_1^7$ is a multiple of $\Lat_1^{4,3}$.

\goodbreak
Our approach also yields a structure theorem for valuations with values in the space of bivariate formal power series $\fps$.
A valuation $\oZb\colon\PZ\to\fps$ is called \Dfn{translatively exponential} if 
\[\oZb(P+v) = \e^v \oZb(P)\]
for every $P\in\PZ$ and $v\in\ZZ^2$, where we again consider $v$ as a linear function on the dual of $\RR^2$. Let $\bVal$ denote the vector space of translatively exponential, $\glz$ equivariant valuations $\oZb\colon\PZ\to\fps$. Note that the functions
\begin{equation}
\label{eq:examples}
    \overline{\Lat}(P)\coloneq  \sum_{v\in P\cap\ZZ^2} \e^v, \qquad \expval(P)\coloneq \int_P \e^v\,\mathrm{d}v
\end{equation}
are both in $\bVal$ (see \cite{Barvinok2008, BeckRobins} for more information and applications).
\goodbreak
For $d\in\ZZ$, we say that $\oZb\colon\PZ\to\fps$ is \Dfn{$d$-dilative} if 
\[\oZb(mP)(x,y) = m^{-d}\oZb(P)(mx,my)\]
for every $m\in\ZZ_{>0}$, $P\in\PZ$, and $x,y\in\RR$.  
For instance, the exponential integral $\expval$ defined in \cref{eq:examples}  is $(-2)$-dilative.
Set
\[
\bVal_d\coloneq \{\oZb\in \bVal\colon \text{ $\oZb$ is $d$-dilative\}}
\]
for $d\in\ZZ$. 
Similar to the translatively polynomial case, dilative valuations decompose the space $\bVal$, and  we determine the dimensions of the vector spaces $\bVal_d$.

\begin{theorem}
\label{thm:exp_val}
For every $\oZb\in\bVal$, there are unique $\oZb_d\in\bVal_d$ with $d\in\ZZ$ such that $\oZb= \sum_{d\geq -2} \oZb_d$.
Moreover,     
    \[ \dim \bVal_d = 
    \begin{cases}
        1&\text{ for }-2\leq d\leq 0\text{ and for }d>0 \text{ odd},\\
       \parts{\frac{d}{2}+1}&\text{ for } d>0 \text{ even},\\
        0&\text{ otherwise.}
    \end{cases}
    \]
\end{theorem}

\noindent
The $r$-summand (that is, the sum of the monomials of degree $r$) of a valuation in  $\bVal_d$ is in $\Val^r_i$ for $d=r-i$.
We will see that all valuations in $\Val_i^r$ for $i>1$ occur as  $(r-i)$-summand of a translatively exponential valuation in $\bVal_{r-i}$.  
For $d\leq 0$ and for $d>0$ odd, the space $\bVal_d$ is spanned by the $d$-dilative part $\overline\Lat_d$ of $\overline\Lat$ (see \cref{sec:lawrence}). For $d>0$ even, we will present in \cref{eq:triang_free} an explicit representation for all valuations in $\bVal_d$. 
Thus, we obtain a triangulation-free description of the valuations in $\Val_i^r$ for $i>1$ and $r-i$ even as summands of valuations in $\bVal_{r-i}$. In \cref{sec:consequences}, we describe some consequences of our results.

\section{Overview of the Proofs}
\label{sec:overview}

In the following sections, we simultaneously prove \cref{thm:tensor_val} and \cref{thm:exp_val}. When the vector spaces in question are zero- or one-dimensional, the results essentially follow from the arguments in \cite{ludwigsilverstein}. We repeat these arguments in \cref{sec:remaining} to keep the proof as self-contained as possible. The majority of the work thus goes in the second and third cases of \cref{thm:tensor_val} and the second case of \cref{thm:exp_val}.

\begin{figure}
    \centering
    \includegraphics[width = .9\textwidth]{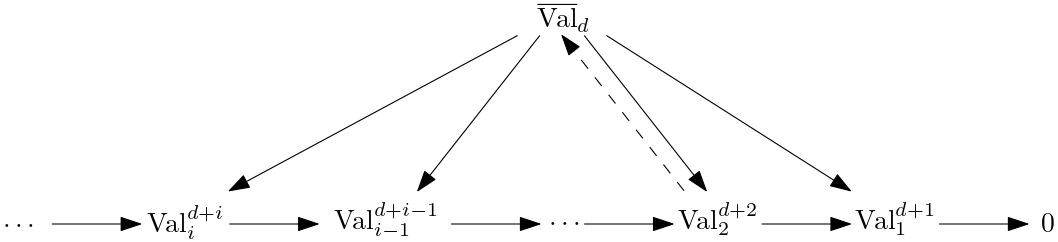}
    \caption{A commutative diagram connecting the spaces of valuations.}
    \label{fig:diagram}
\end{figure}

In the second and third cases of \cref{thm:tensor_val}, the difference $d=r-i$ of rank and homogeneity degree is a positive even integer. If we fix a positive even integer $d$, it turns out that the vector spaces in question are connected by the diagram in \cref{fig:diagram}. In this diagram, the horizontal maps assign to a valuation in $\Val_i^{d+i}$ its first associated valuation, that is, the linear coefficient of $v$ in \cref{eq:translation}. In \cref{sec:pre}, we recapitulate how it follows from \cite{ludwigsilverstein} that the horizontal maps are well-defined. The vertical maps in the diagram simply map a valuation $\oZb\in\bVal$ to its $(d+i)$-summand $\oZ^{d+i}$, where $\oZ^{d+i}(P)$ is the $(d+i)$-summand of the power series $\oZb(P)$ for $P\in\PZ$. The well-definedness of this map is the subject of \cref{lemma:hom_dilative}.

The proofs of \cref{thm:tensor_val} and \cref{thm:exp_val} for even $d>0$ can be summarized as follows:
In \cref{sec:1hom}, we show that $\dim\Val_1^{d+1}=\parts{d+1}$ by relating the respective valuations to invariants of a finite group (see \cref{thm:eval_iso}). In \cref{sec:2hom}, we show that the horizontal maps in the diagram in \cref{fig:diagram} are injective (except for $\Val_1^{d+1}\to 0$, see \cref{lemma:familybusiness}) and compute the dimension of the image of $\Val_2^{d+2}\to\Val_1^{d+1}$. From this, we will obtain that $\dim\Val_2^{d+2} = \parts{\tfrac d2 + 1}$ (see \cref{thm:2hom}).

At this point, it remains to show that the dimensions of $\bVal_d$ and $\Val_i^{d+i}$ for $i>2$ are equal to $\dim\Val_2^{d+2}$. Our strategy for this is as follows: First, in \cref{sec:expo}, we show that the diagram in \cref{fig:diagram} commutes. Then, in \cref{sec:lawrence}, we construct an inverse map to $\bVal_d\to\Val_2^{d+2}$, that is, for an arbitrary $\oZ_2^{d+2}\in\Val_2^{d+2}$, we construct a valuation $\oZb\in\bVal_d$ such that the $(d+2)$-summand of $\oZb(P)$ is $\oZ_2^{d+2}(P)$. The construction is based on an extension to rational, possibly unbounded polyhedra. Using Lawrence's application of the Brion identity \cite{lawrence}, we construct a valuation on rational polyhedra into the space of meromorphic functions whose restriction to lattice polygons takes analytic values.

In \cref{sec:proofs}, we complete the proofs of \cref{thm:tensor_val} and \cref{thm:exp_val}. This section should be regarded as a summary of the previous sections.
Finally, we derive two corollaries from our classification in \cref{sec:consequences}.

\section{Preliminaries}
\label{sec:pre}
We collect results and definitions that are used in our arguments.
 
\subsection{Associated valuations} For $r\in\ZZ_{\ge0}$, let $\oZ\colon\PZ \to \ring_r$ be a trans\-latively polynomial valuation, that is, for any $P\in\PZ$ and $v\in\ZZ^2$, we have
\[
\oZ(P+v) = \sum_{j=0}^r \oZ^{r-j}(P)\frac{v^j}{j!}
\]
with suitable functions $\oZ^{r-j}\colon\PZ\to\ring_{r-j}$. By comparing the coefficients of $v^j$, we see that each $\oZ^{r-j}$ is uniquely determined and a valuation. Moreover, setting $v=0$, we obtain $\oZ^r = \oZ$. We refer to the valuations $\oZ^{r-j}$ as the \Dfn{associated valuations of $\oZ$}. It was shown in \cite{ludwigsilverstein} that associated valuations are, in the following way,  hereditary.

\begin{lemma}[Proposition~6 in \cite{ludwigsilverstein}]
    \label{lemma:hereditary}
    If $\,\oZ\colon\PZ\to\ring_r$ is a translatively poly\-nomial valuation with associated valuation $\oZ^{r-j}\colon\PZ\to\ring_{r-j}$ for $0\leq j\leq r$, then each $\oZ^{r-j}$ is translatively polynomial and its associated valuations are $\oZ^{r-k}$ for $j\leq k \leq r$.
\end{lemma}

It follows from a result of Khovanski\u{\i} and Pukhlikov~\cite{pukhlikovkhovanskii} that $\oZ$ can be decomposed into \Dfn{homogeneous components}. 

\begin{theorem}[Theorem 14 in \cite{ludwigsilverstein}]
    \label{thm:khovanskii}
    If $\,\oZ\colon\PZ\to\ring_r$ is a translatively polynomial valuation, then there exist $i$-homogeneous valuations $\oZ_i\colon\PZ\to\ring_r$ for $0\leq i \leq r+2$ such that $\oZ = \sum_{i=0}^{r+2} \oZ_i$.
\end{theorem}
The next lemma follows from this decomposition and Propositions 6 and 7 in \cite{ludwigsilverstein}.

\begin{lemma}
    \label{lemma:crazywindices}
    Let $\oZ\colon\PZ\to\ring_r$ be a translatively polynomial valuation. For $0\leq i\leq r+2$,  each homogeneous component $\oZ_i$ is translatively polynomial. The associated functions of $\oZ_i$ are given by the homogeneous decomposition of the associated functions of $\oZ$, that is, $(\oZ_i)^{r-j} = (\oZ^{r-j})_{i-j}$ for $j\leq \min\{i,r\}$ and $(\oZ_i)^{r-j} =0$ otherwise.
\end{lemma}

We will use the abbreviation $\oZ^{r-j}_{i-j}\coloneq  (\oZ^{r-j})_{i-j}$.

\subsection{Invariant and equivariant valuations.}
If the valuation $\oZ$ in \cref{lemma:crazywindices} is $\GL_2(\ZZ)$ invariant, then its associated valuations inherit the $\GL_2(\ZZ)$ invariance. For a subgroup $G$ of $\GL_2(\RR)$, let $\ring^G$ be the ring of $G$ invariant polynomials. It is easy to see that
\begin{equation}\label{eq:gl2zinvariants}
    \ring^{\GL_2(\ZZ)}\cong \RR. 
\end{equation}
As observed in \cite{ludwigsilverstein}, this implies that the constant term in the homogeneous decomposition in \cref{thm:khovanskii} vanishes if $\oZ$ is $\GL_2(\ZZ)$ invariant. 

\begin{lemma}
    \label{lemma:constterm}
    For $r\ge 0$ and $\oZ\in\Val^r$, there exists $c\in\RR$ such that $\oZ_0(P)=c$ for every $P\in\PZ$. If $r>0$, then this constant is zero.
\end{lemma}

It was also observed in \cite{ludwigsilverstein} that this implies the following result.
\begin{lemma}
    \label{lemma:onehom_translation_invariant}
     If $r\ge 2$ and $\oZ\in\Val_1^r$, then $\oZ$ is translation invariant. 
\end{lemma}

A valuation $\oZ\colon \PZ\to\ring_r$ is called \emph{odd} if it is equivariant with respect to reflection at the origin, that is, we have $\oZ(-P) = -\oZ(P)$. Note that since the reflection at the origin is in $\glz$, we have $\oZ(-P) = (-1)^r\oZ(P)$ for any $\oZ\in\Val^r$. So $\oZ\in\Val^r$ is odd if and only if $r$ is an odd number.

In the following, we will also deal with valuations taking values in the spaces of formal power series or meromorphic functions on $(\RR^2)^\star$. For such valuations,  the term $\glz$ equivariance will also refer to the action $\phi f = f\circ \phi^\star$, where $\phi\in\glz$ and $f$ is a (formal) function.

\subsection{Planar lattice triangulations}

A \Dfn{decomposition} of a two-dimensional lattice polygon $P$ is a set of lattice polygons $\triangulation$ such that
the following properties hold:
\begin{enumerate}
    \item $P = \bigcup_{T\in\triangulation} T$.
    \item For all $S,T\in\triangulation$, the set $S\cap T$ is a (possibly empty) face of $S$ and $T$.
\end{enumerate}
If all $S\in\triangulation$ are triangles, we call $\triangulation$ a \Dfn{triangulation}.

Recall that a triangle $S=\conv\{v_0,v_1,v_2\}$ with $v_i\in\ZZ^2$ is called \Dfn{unimodular} if it can be expressed as $S = \phi \sttriangle + v$, where $\sttriangle=\conv\{0,e_1,e_2\}$ and $\phi\in\glz$, $v\in\ZZ^2$. Here, $\conv$ stands for convex hull.
The triangle $S$ is unimodular if and only if the edge vectors $v_i-v_0$ for $i\in\{1,2\}$ generate $\ZZ^2$. Note that $v_0$ might be any of the three vertices of $S$ in this condition. 
In fact, it follows that $S$ is unimodular if and only if $\vol(S) = \tfrac 12$, where $\vol$ stands for (two-dimensional) volume. Note that one has $\vol(S) > \tfrac 12$ for a non-unimodular lattice triangle since its edge vectors generate a proper two-dimensional sublattice of $\ZZ^2$.
A triangulation is called unimodular if all of its simplices are unimodular.

Moreover, a lattice triangle is unimodular if and only if it is \emph{empty}; that is, it contains no lattice points other than the vertices. It follows that every lattice polygon $P\in\PZ$ admits an unimodular triangulation.

For triangulations $\triangulation$ and $\triangulation'$ of a lattice polygon $P$ one says that $\triangulation$ and $\triangulation'$ are \emph{connected by a flip}, if $\triangulation\setminus\triangulation' = \{S_1,S_2\}$ and $\triangulation'\setminus\triangulation = \{S_1',S_2'\}$ and, moreover, $S_1\sqcup S_2 = S_1'\sqcup S_2'$ is a convex quadrilateral. In this case, $\{S_1, S_2\}$ and $\{S_1', S_2'\}$ correspond to the two possible triangulations of the quadrilateral.

\begin{theorem}[Lawson \cite{lawson}]
\label{thm:lawson}
    Let $P\in\PZ$ and let $\triangulation$ and $\triangulation'$ be triangulations of $P$. Then there exists a finite sequence of triangulations \(\triangulation = \triangulation_0,\triangulation_1,\dots,\triangulation_m = \triangulation'\) such that $\triangulation_{i-1}$ and $\triangulation_i$ are connected by a flip for all $1\leq i\leq m$.
\end{theorem}

\noindent
Lawson's theorem holds in more generality for triangulations of arbitrary (non-lattice) point sets.
\goodbreak

A valuation $\oZ$ on $\PZ$ is called \emph{simple} if it vanishes on points and segments.
Using the valuation property inductively, it is easy to see that valuations on $\PZ$ satisfy the inclusion-exclusion principle. In particular if $\oZ$ is a simple valuation and $\triangulation$ is a triangulation of a polygon $P\in\PZ$, we have
\begin{equation}
\label{eq:mcmullen}
    \oZ(P) = \sum_{S\in\triangulation} \oZ(S).
\end{equation}
The inclusion-exclusion principle also holds for valuations on lattice polytopes in higher dimensions, but its proof is no longer trivial (see \cite{McMullen09}).

The following corollary of \cref{thm:lawson} helps to construct a valuation ``from scratch''. Let $\rho_i$ denote the reflection on the $i$th coordinate axis and $\US$ the set of unimodular lattice triangles.

\begin{corollary}
    \label{cor:unimodular_flips}
    If $\oz \colon \US\to \ring$ is $\GL_2(\ZZ)$ equivariant and such that \begin{equation}
    \label{eq:well-def}
    \oz(v +\sttriangle) + \oz(v+e_1+e_2 - \sttriangle) = \oz(v+e_1 + \rho_2 \sttriangle) + \oz(v+e_2+\rho_1\sttriangle)
    \end{equation}
    for every $v\in\ZZ^2$, then there is a unique simple, $\GL_2(\ZZ)$ equivariant valuation $\oZ\colon\PZ\to\ring$ with $\oZ(T) = \oz(T)$ for every $T\in\US$. Conversely, if $\oZ\colon\PZ\to\ring$ is a simple, $\GL_2(\ZZ)$ equivariant valuation, then its restriction $\oz$ to $\US$ satisfies \cref{eq:well-def}. 
\end{corollary}

\begin{proof}
    Note that 
    \begin{equation}
    \label{eq:flip}
        \sttriangle\sqcup(e_1+e_2-\sttriangle) = (e_1 + \rho_2 \sttriangle)\sqcup (e_2+\rho_1\sttriangle) = [0,1]^2.
    \end{equation}
    Let $\oz$ be as in the statement of the corollary. For $P\in\PZ$ two-dimensional, we define $\oZ(P;\triangulation) \coloneq  \sum_{T\in \triangulation} \oz(T)$, where $\triangulation$ is a unimodular triangulation of $P$. We claim that $\oZ(P;\triangulation)$ is independent of $\triangulation$. By \cref{thm:lawson}, it is enough to consider a triangulation $\triangulation'$, which is connected to $\triangulation$ by a flip. Let $Q$ be the convex quadrilateral in which $\triangulation$ and $\triangulation'$ differ. It is easy to see that $Q = v+ \phi [0,1]^2$ for some $\phi\in\glz$ and $v\in\ZZ^2$.  
    Since there are exactly two lattice triangulations of $[0,1]^2$, the ones given in \cref{eq:flip}, it follows that, say, $\triangulation$ contains $v + \phi(e_1+e_2-\sttriangle)$ and $v+\phi\sttriangle$. Consequently, $\triangulation'$ contains $v + \phi(e_1 + \rho_2 \sttriangle)$ and $v+\phi(e_2+\rho_1\sttriangle)$. Since $\oz$ is $\GL_2(\ZZ)$ equivariant, we obtain $\oz(v+\phi T) = \oz(\phi^{-1}v +T)$ for any unimodular triangle $T$. Thus, $\oZ(P;\triangulation) = \oZ(P;\triangulation')$ follows from \cref{eq:well-def}, and we can define $\oZ(P) \coloneq  \sum_{T\in \triangulation} \oz(T)$ for $P\in\PZ$ two-dimensional. Setting $\oZ(P)\coloneq 0$ for $\dim P < 2$, we obtain a simple, $\GL_2(\ZZ)$ equivariant valuation with $\oZ(T) = \oz(T)$ for every $T\in\US$. The uniqueness of $\oZ$ follows from \cref{eq:mcmullen}.
    The converse statement is immediate since both sides of \cref{eq:well-def} are equal to $\oZ(v+[0,1]^2)$.
\end{proof}

The following lemma will be used frequently in the following sections.

\begin{lemma}
    \label{lemma:simple_eval}
    Let $\oZ\in\Val^r$ be simple with associated valuations $\oZ^{r-j}$ for $0\leq j\leq r$. 
    If $P\in\PZ$ is two-dimensional with unimodular triangulation 
    $\triangulation = \{ \phi_i \sttriangle + v_i: \phi_i\in\glz,\, v_i\in\ZZ^2,\, 1\leq i \leq m\}$,
    then
    \begin{equation}
    \label{eq:simple_eval}
    \oZ(P) = \sum_{i=1}^m \sum_{j=0}^r  \big(\oZ^{r-j}(\sttriangle)\circ \phi_i^\star\big) \frac{v_i^j}{j!}.
\end{equation}
    In particular, $\oZ$ is uniquely determined by the values $\oZ^{r-j}(\sttriangle)$ for $1\leq j\leq r$.
\end{lemma}

\begin{proof}
    By using the equivariance properties of $\oZ$ to \cref{eq:mcmullen}, we immediately obtain \cref{eq:simple_eval}. For the last statement, we recall that every two-dimensional lattice polygon has a unimodular triangulation, so $\oZ(P)$ can always be evaluated via \cref{eq:simple_eval}.
\end{proof}

\subsection{Integer partitions}

In Theorem \ref{thm:tensor_val}, we defined
\begin{equation}\label{eq:part}
    \parts{r} \coloneq  |\{(k,\ell)\in\ZZ_{\geq 0}^2 \colon 2k+3\ell=r\}|,
\end{equation}
the number of integer partitions of $r$ into parts $2$ and $3$. In explicit terms, one has
\begin{equation}
    \label{eq:p23explicit}
    \parts{r} = \left\lfloor \frac{r+2}{2}\right\rfloor - \left\lfloor \frac{r+2}{3}\right\rfloor.
\end{equation}
The sequence is \href{https://oeis.org/A103221}{A103221} in the {Online Encyclopedia of Integer Sequences}.

\section{One-homogeneous valuations of odd rank}
\label{sec:1hom}

In this section, we prove \cref{thm:tensor_val} for $i=1$ and $r>1$ odd.
We require the following result. 

\begin{lemma}
\label{lemma:simple}
    Let $r-i>0$ be even and $i\ge 1$. If $\,\oZ\in\Val_i^r$, then $ \oZ$ is simple.
\end{lemma}

\begin{proof}
    First, let $i=1$. It follows from the assumption that $r$ is odd and from \cref{lemma:onehom_translation_invariant}, that $\oZ$ is translation invariant. Let $Q\in\PZ$ be a segment or a point. In either case, a vector $v\in\ZZ^2$ exists, such as $-Q = Q+v$. Since $r$ is odd and $\oZ$ is $\glz$ equivariant, the valuation $\oZ$ is also odd. 
    Hence, by translation invariance,
    \[
    -\oZ(Q) = \oZ(-Q)=\oZ(Q+v) = \oZ(Q)
    \]
    and, thus, $\oZ(Q)=0$.

\goodbreak
    Now, let $i>1$. By induction and \cref{lemma:constterm}, the associated functions $\oZ_{i-j}^{r-j}$ of $\oZ$ are simple for $0<j<i$. This implies $\oZ(Q+v)=\oZ(Q)$ for any $v\in\ZZ^2$ and any lower-dimensional $Q\in\PZ$.

    It follows from \eqref{eq:gl2zinvariants} that $\oZ(\{0\})=0$. For any $p\in\ZZ^2$ with $p\ne 0$,  the translation invariance on lower-dimensional polygons implies that $\oZ(\{p\}) = \oZ(p+\{0\})=0$. So, let $S\in\PZ$ be a segment. For $m\ge 1$, we can decompose the dilate $mS$ into $m$ translates $S_i$ with $1\leq i\leq m$ of $S$ that intersect at lattice points. We just saw that our valuation vanishes on points. Hence, we obtain from the homogeneity of $\oZ$ that
    \[
    m^i\oZ(S) = \oZ(mS) = \sum_{i=1}^m \oZ(S_i) = m\oZ(S).
    \]
    Since $i>1$, we conclude that $\oZ(S)=0$.
\end{proof}

Let $\oZ\in\Val_1^r$ with $r>1$ odd. As we saw earlier, $ \oZ$ is translation invariant by \cref{lemma:onehom_translation_invariant}, simple by \cref{lemma:simple}, and odd. The converse statement is also true.

\begin{lemma}
    \label{lemma:translation_invariant_onehom}
    If $r$ is odd and $ \oZ\colon\PZ\to\ring_r$ is a simple, translation invariant, $\GL_2(\ZZ)$ equivariant valuation, then $ \oZ$ is one-homogeneous, that is, $\oZ\in \Val_1^r$.
\end{lemma}

\begin{proof}
    Since the valuation $\oZ$ is translation invariant, all its associated valuations $\oZ^{r-j}\colon\PZ\to\ring_{r-j}$ vanish except for $\oZ^r=\oZ$. Thus, \cref{lemma:simple_eval} implies that it is enough to check $\oZ(m\sttriangle) = m\oZ(\sttriangle)$ for $m\in\ZZ_{>1}$.

\goodbreak
    We observe that \[
    m\sttriangle = \Big(\bigsqcup_{v\in(m-2)\sttriangle\cap\ZZ^2} (v+[0,1]^2)\Big) 
    \sqcup \Big( \bigsqcup_{a,b\in\ZZ_{\geq 0},~a+b=m-1} (ae_1 +be_2 + \sttriangle) \Big).
    \]
    Since $\oZ$ is odd and translation invariant and $[0,1]^2$ is centrally symmetric, it follows that $\oZ(v+[0,1]^2) = 0$ for all $v\in\ZZ^2$. Thus, by \cref{eq:mcmullen},
    \[
    \oZ(m\sttriangle) = \sum_{a,b\in\ZZ_{\geq 0},~a+b=m-1} \oZ(ae_1 + be_2 +\sttriangle) = m\oZ(\sttriangle),
    \]
    where we used again that $\oZ$ is simple and translation invariant.
\end{proof}

\goodbreak
Let $\group$ denote the subgroup of unimodular transformations $\phi\in\glz$ for which there exists $v\in\ZZ^2$ such that $\sttriangle = \phi\sttriangle + v$.

\begin{theorem}
\label{thm:eval_iso}
    Let $r>1$ be odd. The evaluation map
    \[
    \eval^r \colon \Val_1^r\to \ring_r^\group,\quad \eval^r(\oZ) \coloneq  \oZ(\sttriangle)
    \]
    is a well-defined linear isomorphism. The valuation $\oZ_f: =(\eval^r)^{-1}(f)$, obtained by the inverse map, is simple for $f\in\ring_r^\group$ and given by
    \begin{equation}
          \label{eq:adhoc_definition}
            \oZ_f(P) = \sum_{i=1}^m f\circ \phi_i^\star.
    \end{equation}
    for any two-dimensional $P\in\PZ$ 
     with unimodular triangulation 
    $\triangulation = \{ \phi_i \sttriangle + v_i: \phi_i\in\glz,\, v_i\in\ZZ^2,\, 1\leq i \leq m\}$.
\end{theorem}

\begin{proof}
    First, we show that $\eval^r(\oZ)\in\ring_r^\group$ for $\oZ\in\Val_1^r$. If $\phi\in \group$, then there exists $v\in\ZZ^2$ such that $\sttriangle=\phi\sttriangle + v$. Since $\oZ$ is translation invariant by \cref{lemma:onehom_translation_invariant} and $\GL_2(\ZZ)$ equivariant, it follows that
    \(
    \oZ(\sttriangle) = \oZ(\sttriangle) \circ\phi^\star
    \)
    for $\phi\in\group$,
    which shows that $\oZ(\sttriangle)\in\ring_r^\group$ as required. 
    
    The linearity of $\Theta^r$ is obvious. The injectivity follows directly from \cref{lemma:simple_eval}. To show that $\eval^r$ is surjective, we
    show that for any $f\in\ring_r^\group$ the expression \cref{eq:adhoc_definition} can serve as a definition of a simple valuation $\oZ_f\in\Val_1^r$, which is then necessarily $(\Theta^r)^{-1}(f)$.

    To this end, 
    consider a unimodular triangle $T\in\PZ$ represented in two ways as an affine image of the standard triangle $\sttriangle$, that is, $T = \phi\sttriangle + v = \phi'\sttriangle+v'$ for certain $\phi,\phi'\in\glz$ and $v,v'\in\ZZ^2$.
    Rearranging yields \[
    \sttriangle = \phi^{-1}\phi'\,\sttriangle +\phi^{-1}(v'-v).
    \]
    Hence, $\phi^{-1}\phi'\in \group$. Since $f$ is $\group$ invariant, it follows that $f\circ (\phi')^\star = f\circ\phi^\star$, that is, $\oz(T)\coloneq  f\circ\phi^\star$ is a well-defined translation invariant, $\GL_2(\ZZ)$ equivariant function on the set of unimodular triangles.

    Next, we note that the triangles $\rho_1\sttriangle$ and $\rho_2\sttriangle$ occurring in \cref{eq:well-def} are opposite in the sense that $\rho_1\sttriangle = -\rho_2\sttriangle$. Thus, since the degree of $f$ is odd, we have $\oz(-\sttriangle) = -\oz(\sttriangle)$ and $\oz(\rho_1\sttriangle) = -\oz(\rho_2\sttriangle)$. Since $\oz$ is translation invariant, \cref{eq:well-def} is satisfied. Hence, \cref{cor:unimodular_flips} yields a simple $\glz$ equivariant valuation $\oZ\colon\PZ\to\ring_r $ with $\oZ(T) = \oz(T)$ for all unimodular triangles $T$. Since $\oz$ is translation invariant, it follows from \cref{eq:mcmullen} that $\oZ$ is translation invariant as well. In particular, $\oZ\in\Val^r$ and \cref{lemma:simple_eval} shows that $\oZ(P) = \oZ_f(P)$ as given by \cref{eq:adhoc_definition}.

    Finally, \cref{lemma:translation_invariant_onehom} yields that $\oZ_f$ is 1-homogeneous, that is, we have $\oZ_f\in\Val_1^r$ as desired. 
\end{proof}

For $r=1$, the valuations in $\Val_1^1$ are not necessarily translation invariant. Thus, $\oZ(\sttriangle)$ is not $\,\group$ invariant in this case. For $r>1$ even, the valuations in $\Val_1^r$ are translation invariant by \cref{lemma:onehom_translation_invariant}, and so the evaluation map $\eval^r$ is a well-defined linear map to $\ring_r^\group$. However, it is, in general, not an isomorphism.

To compute the dimension of $\Val_1^r$ for $r>1$ odd, it is now enough to compute the Molien series of the invariant ring $\ring^\group$. To this end, we
observe that $\group$ has order 6. Any unimodular transformation $\phi\in\GL_2(\ZZ)$ maps $\sttriangle$ to a triangle with a vertex at the origin. So if $\phi\sttriangle=\sttriangle+v$, we have $v\in\{0,-e_1,-e_2\}$. The two non-zero vertices of $\sttriangle$ must then be mapped by $\phi$ to the non-zero vertices of $\sttriangle+v$. Hence, we have
\begin{equation*}
\label{eq:gens}
\group = \left\langle
\begin{pmatrix}
    0 & 1\\
    1 & 0
\end{pmatrix},
\begin{pmatrix}
    -1 & -1\\
   \phantom{-}0 & \phantom{-}1
\end{pmatrix}
\right\rangle,
\end{equation*}
since the group on the right-hand side has order 6 and is contained in $\group$. From this, we see that $\group$ is a finite group generated by \emph{pseudo-reflections}, that is, linear transformations that fix a space of codimension~1. The Chevalley--Shephard--Todd Theorem (see, for example, \cite[Theorem 7.1.4]{MR1869812}) now implies that the ring of $\group$ invariant polynomials, $\ring^\group$, is finitely generated and that its generators are algebraically independent. Molien's formula (for example, \cite[Theorem 3.1.3]{MR1869812}) yields that the dimension of $\ring_r^\group$ is $\parts{r}$ (defined in \cref{eq:part}).  In particular, $\ring^\group$ is generated by suitable polynomials of degree $2$ and $3$.

\goodbreak
\begin{lemma}
    \label{lemma:chevalley}
    The algebraically independent polynomials 
      $\ptwo(x,y) \coloneq  x^2 -xy+y^2$ and $\pthree(x,y) \coloneq x^3-\tfrac 32(x^2y+xy^2)+y^3$
      generate the invariant ring $\ring^\group$.
\end{lemma}

\begin{proof}
    It is enough to apply the Reynolds operator $\frac{1}{\lvert\group\rvert}\sum_{\phi\in\group} \phi$, which is a projection onto $\ring^\group$, to the polynomials $x^2$ and $x^3$ and thus compute the generators $\ptwo$ and $\pthree$.
\end{proof}

\goodbreak
Combining \cref{lemma:chevalley} with \cref{thm:eval_iso} gives the second case of \cref{thm:tensor_val}.

\begin{corollary}
    \label{cor:onehom_dimension}
    For $r>1$ odd, we have $\dim\Val_1^r = \parts{r}$.
\end{corollary}

Note that $r=9$ is the first odd number such that $\parts{r} > 1$.
Since $\Lat^9(m\sttriangle)$ is a polynomial of degree 11 in $m$ with coefficients $\Lat_i^9(\sttriangle)$, we can compute $\Lat_1^9(\sttriangle)$ via interpolation. 
Using {\tt SageMath} \cite{sagemath}, we obtain $\Lat_1^9 = \oZ_f$ (defined in \cref{eq:adhoc_definition}), where 
\[
  f = \frac{1}{990\cdot 9!}\left( 4\pthree^3 - 79\ptwo^3\pthree\right).
\]
The valuation constructed in \cite[Section 8.2]{ludwigsilverstein} is a unimodular valuation that is not a linear combination of Ehrhart tensor coefficients. It corresponds (up to a constant factor) to $\oZ_g$ with $g= \pthree^3$.

\section{Two-homogeneous valuations of even rank}\label{sec:two-hom}
\label{sec:2hom}

Throughout the section, we assume that $r>1$ is odd and consider a valuation $\oZ_1^r\in\Val_1^r$. When does there exist  
$\oZ_2^{r+1}\in\Val^{r+1}$ such that 
\begin{equation}
\label{eq:2trans}
    \oZ_2^{r+1}(P+v) = \oZ_2^{r+1}(P) + \oZ_1^r(P)v
\end{equation}
for every $v\in\ZZ^2$ and $P\in\PZ$?

To answer this question, we first prove the following lemma on associated valuations, which we will also use in the classification for higher homogeneity degrees.

\begin{lemma}
\label{lemma:familybusiness}
   Let $\oY, \oZ\in\Val_i^r$ with $r>i>1$ and $r-i$ even. If the associated valuations satisfy $\,\oY^{r-j}=\oZ^{r-j}$ for all $0<j<i$, then $\oY=\oZ$.
\end{lemma}

\goodbreak
\begin{proof}[Proof of \cref{lemma:familybusiness}]
    Let $\oX \coloneq  \oZ-\oY$. By assumption, $\oX$ is translation invariant, as well as $\GL_2(\ZZ)$ equivariant and $i$-homogeneous, and by \cref{lemma:simple}, it is simple. If $i=2$, then \cite[Proposition 23]{ludwigsilverstein} yields $\oX=0$ and thus the claim. If $i>2$, we have $X=0$ by \cite[Theorem 5]{McMullen77}.  
\end{proof}

By \cref{lemma:familybusiness}, there is at most one solution $\oZ_2^{r+1}$ of \eqref{eq:2trans}. We will now determine for which $\oZ_1^r$ a solution $\oZ_2^{r+1}$ exists. We write $\sqgroup$ for the group generated by the reflections at the coordinate axes and permutation of the coordinates, that is, 
  \[
    \sqgroup=\left\langle \begin{pmatrix}
        0 & 1\\
        1 & 0
    \end{pmatrix},\,\begin{pmatrix}
        1 & \phantom{-}0\\
        0 & -1
    \end{pmatrix}\right\rangle.
    \]
Note that $\sqgroup$ is a finite subgroup of $\glz$.

\begin{proposition}
\label{prop:2hom_equiv}
    Let $r>1$ be odd and $\oZ_1^r\in\Val_1^r$. The following are equivalent:
    \begin{enumerate}
        \item There exists $\oZ_2^{r+1}\in\Val_2^{r+1}$ such that \cref{eq:2trans} holds for every $P\in\PZ$ and $v\in\ZZ^2$.
        \item\label{it:invariance} There exists  $h\in\ring_{r+1}^\group$ such that $h+(x+y)\oZ_1^r(\sttriangle)\in\ring_{r+1}^\sqgroup$.
    \end{enumerate}
    Moreover, the polynomial $h$ is uniquely determined by $\oZ_1^r(\sttriangle)$ if the conditions \emph{(1)} and \emph{(2)} are satisfied.
\end{proposition}

\begin{proof}
    First, we observe that condition (2) is by linearity equivalent to
    \begin{itemize}
        \item[(2')] There exists  $h\in\ring_{r+1}^\group$ such that $2h-(\tfrac{x}{3}+\tfrac y3)\oZ_1^r(\sttriangle)\in\ring_{r+1}^\sqgroup$.
    \end{itemize}
    In the proof, we will work with (2') rather than (2).

    We start by showing that (1) implies (2').
    By \cref{lemma:simple}, the valuation $\oZ_2^{r+1}$ is simple. So, it is determined by  
    $$f_1 \coloneq  \oZ_1^r(\sttriangle) \text{ and }f_2\coloneq \oZ_2^{r+1}(\sttriangle)$$ 
    by \cref{lemma:simple_eval}. Consider a unimodular transformation $\phi\in\group$ such that $\phi\sttriangle + v = \sttriangle$ for some $v\in\ZZ^2$.
    Then, by \cref{eq:2trans} and since $f_1\in\ring_r^\group$,
    \begin{equation}
    \label{eq:f2vsg}
        f_2 = \oZ_2^{r+1}(\phi \sttriangle + v) = f_2\circ\phi^\star + f_1 v.
    \end{equation}
    Suppose $\widetilde{f}_2$ is another solution of \eqref{eq:f2vsg}. Then $f_2-\widetilde{f}_2\in\ring_{r+1}^\group$. Next, we make the guess that $(\tfrac{x}{3}+\tfrac{y}{3})f_1$ is a particular solution to \cref{eq:f2vsg}. To see that the guess is correct, note that $\centroid(\sttriangle)\coloneq \tfrac{x}{3} + \tfrac{y}{3}$ corresponds to the centroid of the standard triangle. 
    \[
    \big(\centroid(\sttriangle)f_1\big)\circ \phi^\star + vf_1 = \centroid(\phi \sttriangle)(f_1\circ\phi^\star) +vf_1 = \centroid(\phi\sttriangle +v) f_1=\centroid(\sttriangle)f_1.
    \]
    So $(\tfrac{x}{3}+\tfrac{y}{3})f_1$ is indeed a solution to \cref{eq:f2vsg} and $f_2 = h+(\tfrac{x}{3}+\tfrac{y}{3})f_1$ for a suitable $h\in\ring_{r+1}^\group$.

    Now consider the square $[0,1]^2$. This is a centrally symmetric polygon, so we have $\oZ_1^r([0,1]^2)=0$ because $\oZ_1^r$ is odd and translation invariant.
    From \cref{eq:2trans}, it follows that $\oZ_2^{r+1}([0,1]^2+v)=\oZ_2^{r+1}([0,1]^2)$ for any $v\in\ZZ^2$, which implies that 
    \begin{equation}\label{eq_contained}
    \oZ_2^{r+1}([0,1]^2)\in\ring_{r+1}^\sqgroup,
    \end{equation}
    since $\phi [0,1]^2$ is a translate of $[0,1]^2$ for $\phi\in\sqgroup$. Using the triangulation  $[0,1]^2 =   \sttriangle\sqcup (-\sttriangle+e_1+e_2)$, we obtain
    \begin{align*}
    \oZ_2^{r+1}([0,1]^2) &= \oZ_2^{r+1}(\sttriangle) + \oZ_2^{r+1}(-\sttriangle+e_1+e_2)\\
    &= f_2(x,y) + f_2(-x,-y) +(x+y)f_1(-x,-y)\\
    &= 2f_2(x,y)-(x+y)f_1(x,y) \\
    &= 2h-(\tfrac{x}{3}+\tfrac{y}{3})f_1.
    \end{align*}
Combined with \eqref{eq_contained}, this proves the implication.

Second, we show that (2') implies (1).    
We define $f_2\coloneq (\tfrac{x}{3}+\tfrac{y}{3})f_1 + h$. For 
a unimodular triangle $T=\phi \sttriangle +v$ with $\phi\in\glz$ and $v\in\ZZ^2$, we set
\begin{equation*}
    \label{eq:def_z}
    \oz(T) \coloneq  f_2\circ\phi^\star + v\cdot f_1\circ\phi^\star.
\end{equation*}
We claim that this is independent of the representation of $T$ as an affine image of $\sttriangle$. As in the proof of \cref{thm:eval_iso}, consider
 $\phi \sttriangle +v = \psi\sttriangle+w$ for suitable $\phi,\psi\in\GL_2(\ZZ)$ and $v,w\in\ZZ^2$. Then, $\sttriangle = \phi^{-1} \psi\sttriangle + \phi^{-1}(w-v)$. This means $\phi^{-1}\psi\in\group$. Set $f_0 \coloneq  (\tfrac x3 + \tfrac y3)f_1$. We already saw in the previous step that $f_0$ satisfies \[
f_0 = f_0\circ (\phi^{-1}\psi)^\star + \phi^{-1}(w-v)\cdot f_1 = f_0\circ\psi^\star\circ\phi^{-\star} + \phi^{-1}(w-v)\cdot f_1.
\]
Adding $h$ to the equation, applying $\phi$ and rearranging  gives
\[\begin{split}
f_2\circ\phi^\star + v\cdot f_1\circ\phi^\star &= h\circ\phi^\star + f_0\circ\phi^\star + v\cdot f_1\circ\phi^\star\\
&= h\circ \phi^\star + f_0\circ\psi^\star + w\cdot(f_1\circ \phi^\star)\\
&= h\circ (\phi^{-1}\psi)\circ\phi^\star + f_0\circ\psi^\star + w\cdot(f_1\circ (\phi^{-1}\psi)\circ \phi^\star)\\
&=h\circ \psi^\star + f_0\circ\psi^\star + w\cdot(f_1\circ \psi^\star)\\
&=f_2\circ\psi^\star + w\cdot(f_1\circ\psi^\star)
\end{split}\]
Thus, $\oz(T)$ is, in fact, well-defined. Moreover, it is easy to verify that $\oz$ is $\GL_2(\ZZ)$ equivariant.

For $v\in\ZZ^2$, consider the triangles $T_1 = v+\sttriangle$ and $T_2 = v+e_1+e_2-\sttriangle$, as well as $T_3=v+e_1+\rho_2\sttriangle$ and $T_4=v+e_2+\rho_1\sttriangle$, where again $\rho_i\in\sqgroup$ denotes the reflection at the $i$th coordinate axis. Using that $r$ is odd, we have
\[
\begin{split}
    \oz(T_1) + \oz(T_2) &= f_2 + vf_1 + f_2 - (x+y)f_1 - vf_1\\
    &= 2(h+(\tfrac x3+ \tfrac y3)f_1) - (x+y)f_1\\
    &= 2h - (\tfrac x3 + \tfrac y3) f_1.
\end{split}
\]
Using that $r$ is odd, we also obtain
\[
\begin{split}
    \oz(T_3) + \oz(T_4) &= f_2(-x,y) + (x+v)f_1(-x,y) + f_2(x,-y) + (y+v)f_1(x,-y)\\
    &= 2f_2(-x,y) + (x-y)f_1(-x,y)\\
    &= 2h(-x,y) - (-\tfrac x3 +\tfrac y3) f_1(-x,y)\\
    &= 2h - (\tfrac x3 + \tfrac y3) f_1.
\end{split}
\]
We used the $\sqgroup$ invariance of $2h - (\tfrac x3 + \tfrac y3) f_1$ to obtain the last line. So the function $\oz$ satisfies \cref{eq:well-def}. By \cref{cor:unimodular_flips}, we find a simple, $\GL_2(\ZZ)$ equivariant valuation $\oZ\colon \PZ \to \ring_{r+1}$ with $\oZ(T) = \oz(T)$ for all unimodular triangles $T$. In fact, $\oZ$ satisfies \cref{eq:2trans}: Since $\oZ$ and $\oZ_1^r$ are simple, it suffices to check \cref{eq:2trans} for a unimodular triangle $T = \phi\sttriangle + w$ and a translation vector $v\in\ZZ^2$. We have
\[
\begin{split}
\oZ(T+v) &= \oz(T+v) = f_2\circ \phi^\star + (v+w)\cdot f_1\circ\phi^\star\\
&= \oz(T) + v\cdot \oZ_1^r(\sttriangle) \circ\phi^\star = \oZ(T) + v\oZ_1^r(T),
\end{split}
\]
where in the last step we used that $\oZ_1^r$ is translation invariant by \cref{lemma:onehom_translation_invariant}. It remains to show that $\oZ$ is 2-homogeneous.
 First, we show that $\oZ(2^k T) = 2^{2k}\oZ(T)$ holds all unimodular simplices $T$ and for all $k\in\ZZ_{\geq 0}$ by induction on $k$. For $k=0$, there is nothing to prove, so let $k>0$. Multiplying the triangulation     
 \begin{equation*}
 \label{eq:standard_triangulation}
       2\sttriangle = \sttriangle\sqcup (\sttriangle+e_1)\sqcup (\sttriangle+e_2) \sqcup (-\sttriangle+e_1+e_2).
\end{equation*}
with $2^{k-1}$, we obtain a triangulation of $2^k\sttriangle$ into four copies of $\pm 2^{k-1}\sttriangle$:
\begin{equation*}
    \label{eq:kst_triangulation}
    2^k\sttriangle = 2^{k-1}\sttriangle \sqcup 2^{k-1}(\sttriangle + e_1) \sqcup 2^{k-1}(\sttriangle+e_2) \sqcup 2^{k-1}(-\sttriangle+e_1+e_2).
\end{equation*}
Using induction, we find
\begin{align*}
\oZ&(2^k\sttriangle) \\
&= \oZ(2^{k-1}\sttriangle) + \oZ(2^{k-1}(\sttriangle + e_1))+\oZ(2^{k-1}(\sttriangle+e_2)) + \oZ(2^{k-1}(-\sttriangle+e_1+e_2))\\
&= 2^{2k-2}\big( \oZ(\sttriangle) + \oZ(\sttriangle+e_1)
+ \oZ(\sttriangle+e_2) + \oZ(-\sttriangle+e_1+e_2)\big)\\
&= 2^{2k}\oZ(\sttriangle),
\end{align*}
where we used \cref{eq:2trans} to obtain the last line.
For any $\phi\in\GL_2(\ZZ)$, we then also have $\oZ(2^k\phi\sttriangle) = 2^{2k}\oZ(\phi\sttriangle)$. For an arbitrary unimodular triangle $\phi\sttriangle+v$, it follows from \cref{eq:2trans} that
\[\begin{split}
\oZ(2^k(\phi\sttriangle+v)) &= \oZ(2^k \phi \sttriangle) + \oZ_1^r(2^k\phi\sttriangle) 2^kv \\
&=2^{2k} (\oZ(\phi\sttriangle) + \oZ_1^r(\phi\sttriangle)v) = 2^{2k}\oZ(\phi\sttriangle+v),
\end{split}\]
where we used that $\oZ_1^r$ is one-homogeneous. So we have $\oZ(2T) = 2^{2k}\oZ(T)$ for all unimodular triangles $T$ and all $k\in\ZZ_{\geq 0}$. By \cref{thm:khovanskii},  this means that $\oZ$ is 2-homogeneous on unimodular triangles. It follows from \cref{eq:mcmullen} that $\oZ$ is 2-homogeneous on arbitrary $P\in\PZ$.

Finally, we establish the uniqueness of $h$.
Note that in the above construction we have $$\oZ([0,1]^2) = 2{h}-(\tfrac{x}{3}+\tfrac y3)f_1.$$
Suppose that a second $\sqgroup$ invariant polynomial $\Tilde{h}$ exists such that 
$$2\Tilde{h}-(\tfrac{x}{3}+\tfrac y3)f_1\in\ring_{r+1}^\sqgroup.$$ 
Then there exists a valuation $\oZt\in\Val_2^{r+1}$ that satisfies \cref{eq:2trans} and has $$\oZt([0,1]^2) = 2\Tilde{h}-(\tfrac{x}{3}+\tfrac y3)f_1
\neq \oZ([0,1]^2).$$ 
However, by \cref{lemma:familybusiness}, we have $\oZ=\oZt$, a contradiction.
\end{proof}

\goodbreak
\cref{prop:2hom_equiv} allows us to determine the dimension of $\V_2^{r+1}$ as follows.
By \cref{lemma:familybusiness} and \cref{thm:eval_iso}, the linear map
\[
\V_2^{r+1}\to  \ring_r^\group,\quad \oZ_2^{r+1}\mapsto \oZ_1^r(\sttriangle)
\]
is injective. \cref{prop:2hom_equiv} tells us that its image consists of those  
$f\in\ring_r^\group$ for which we can find 
$h\in\ring_{r+1}^\group$ such that $h+(x+y)f\in\ring_r^\sqgroup$.

We consider the linear maps
\begin{equation*}
\label{eq:rho}
\rho \colon  \ring_r^\group\times \ring_{r+1}^\group \to \ring_{r+1},\quad (f,h)\mapsto h+(x+y)f
\end{equation*}
and 
\begin{equation*}
\label{eq:delta}
\delta \colon \ring_{r+1} \to \ring_{r+1},\quad g\mapsto \tfrac 12 ( g(x,y)-g(x,-y) ).
\end{equation*}
Since the permutation of the coordinates $x$ and $y$ is in both $\group$ and $\sqgroup$, it follows for $f,h\in\ring^\group$ that $h+(x+y)f\in\ring^\sqgroup$ if and only if $h+(x+y)f\in\ker(\delta)$, the kernel of $\delta$. 
We have 
\begin{equation}\label{eq:dimker}
\dim \V_2^{r+1} = \dim\ker(\delta\circ \rho),    
\end{equation} 
since the natural projection $\ker(\delta\circ \rho)\to\ring_r^\group$ is injective by \cref{prop:2hom_equiv}.

\goodbreak
In the remainder of the section, we prove the following result.
\begin{lemma}
\label{lemma:dimim}
    For $r>1$ odd, the dimension of the image of $\delta\circ\rho$ is $\lfloor \tfrac{r+3}{4}\rfloor$.
\end{lemma}

Using \cref{eq:p23explicit}, it can be checked that \(\parts{r} + \parts{r+1}-\big\lfloor \tfrac{r+3}{4}\big\rfloor = \parts{\tfrac{r+1}{2}}\) for $r$ odd.  Since 
\(\parts{r}+\parts{r+1} = \dim(\ring_r^\group \times \ring_{r+1}^\group)\), 
we obtain the following result from \cref{lemma:dimim} and \cref{eq:dimker}.

\begin{proposition}
    \label{thm:2hom}
    For $r>1$ odd,  $\dim\V_2^{r+1} = \parts{\frac{r+1}{2}}$.
\end{proposition}

For the proof of \cref{lemma:dimim}, we use a third matrix group,
\[
\thirdgroup \coloneq  \left\langle \begin{pmatrix}
    0 & 1 \\
    1 & 0 
\end{pmatrix},
\begin{pmatrix}
    -1 & \phantom{-}0\\
    \,\,0 & -1
\end{pmatrix}
\right\rangle,
\]
a subgroup of $\sqgroup$.
The algebraically independent polynomials 
\[\gprd \coloneq  x y \,\,\text{ and }\,\, \gsqr\coloneq x^2+y^2\]
generate the invariant ring  $\RR[x,y]^\thirdgroup$.

\begin{lemma}
    \label{lemma:imrho}
   For $r>1$ odd,  the image of $\rho$ is in $\ring^{\thirdgroup}_{r+1}$.  As a vector space,
  it is spanned by the polynomials
  \[
    \Big\{(\gsqr - \gprd)^k 
      (2\gsqr - 5\gprd)^\ell 
      (\gsqr + 2\gprd)^{\lfloor\frac{\ell+1}{2}\rfloor} 
      : k = \frac{r+1}{2} - \ell - \Big\lfloor\frac{\ell+1}{2}\Big\rfloor,\, 
      0\leq\ell\leq\Big\lfloor\frac{r+1}{3}\Big\rfloor 
    \Big\}.
    \]
\end{lemma}

\begin{proof}
  We observe that
  \begin{align*}
   \ptwo &= \gsqr - \gprd,\\
   (x+y)\pthree &= \frac{1}{2}(x+y)^2(2x^2 - 5xy + 2y^2) \\
                &= \frac{1}{2}(\gsqr + 2\gprd)(2\gsqr - 5\gprd),\\
   \pthree^2 &= \frac{1}{4}(\gsqr + 2\gprd)(2\gsqr - 5\gprd)^2.
  \end{align*}
  Therefore, any element in the image of $\rho$ is a linear combination of terms of the form $(\gsqr - \gprd)^k 
  (2\gsqr - 5\gprd)^\ell 
  (\gsqr + 2\gprd)^{\lfloor\frac{\ell+1}{2}\rfloor}$ such that $2\big(k+\ell+\lfloor\frac{\ell+1}{2}\rfloor\big) = r+1$.
  
  To determine the range of $\ell$, note that $k$ must be non-negative, and therefore 
  \[
    \ell + \left\lfloor\frac{\ell+1}{2}\right\rfloor \leq \frac{r+1}{2}.
  \]
  By distinguishing cases, we see that this is equivalent to
  \[
    \ell\leq\left\lfloor\frac{r+1}{3}\right\rfloor,
  \]
  provided that $r$ is odd.
\end{proof}

\goodbreak
\begin{lemma}
\label{lemma:imdelta}
 For $r>1$ odd, we have
   \[
    \delta(\gprd^i \gsqr^j) =
    \begin{cases}
      \gprd^i\gsqr^j  &\text{for $i$ odd,}\\
      0 & \text{otherwise,}
    \end{cases}
  \]
  for $0\leq i \leq \tfrac{r+1}{2}$ and $j\coloneq \tfrac{r+1}{2}-i$.
 In particular, $\delta(\ring_{r+1}^\thirdgroup)\subseteq\ring^{\thirdgroup}_{r+1}$.
\end{lemma}
\begin{proof}
  It suffices to observe that $\gsqr$ is even in $y$ and $\gprd$ is odd in $y$.
\end{proof}

\begin{proof}[Proof of \cref{lemma:dimim}]
For $r>1$ odd, let $d(r)$ be the dimension of the image of $\delta\circ\rho$. 
Set
\[\gg_\ell \coloneq  (\gsqr - \gprd)^k 
 (2\gsqr - 5\gprd)^\ell 
 (\gsqr + 2\gprd)^{\lfloor\frac{\ell+1}{2}\rfloor},\]
with $k \coloneq  \tfrac{r+1}{2} - \ell - \lfloor\frac{\ell+1}{2}\rfloor$ for $0\leq\ell\leq\lfloor\frac{r+1}{3}\rfloor$.  
It follows from \cref{lemma:imrho} that $d(r)$ is the number of linearly independent polynomials in 
\[
\Big\{\delta(\gg_\ell): 0\leq\ell\leq\Big\lfloor\frac{r+1}{3}\Big\rfloor\Big\}.
\]
Let $R$ be the matrix of coefficients of monomials of these polynomials, that is,  $R_{\ell, i}$ is the coefficient of $\gprd^i\gsqr^{\frac{r+1}{2}-i}$ in $\delta(\gg_\ell)$.  Since  $\delta(\gg_\ell)$ has no monomials of even degree in $\gprd$ by \cref{lemma:imdelta},   every second column of $R$ vanishes.  We will show that the remaining \[\Big\lfloor\frac{r+3}{4}\Big\rfloor = \Big\lfloor\frac{\frac{r+1}{2} + 1}{2}\Big\rfloor\] columns of $R$ are linearly independent.  This determines the dimension $d(r)$.

By \cref{lemma:imdelta}, the entries of $R$ are integers. 
Therefore, it suffices to show that the remaining columns of $R$ are linearly independent modulo~$2$. This is easier, since we have \[\gg_\ell = (\gsqr + \gprd)^k \gprd^\ell\gsqr^{\left\lfloor\frac{\ell+1}{2}\right\rfloor}\] in $\mathbb Z_2[\gsqr, \gprd]$.  By the binomial theorem, we obtain that $R_{\ell, 2i+1}$ equals, modulo~$2$, 
  \[
    \binom{k}{2i+1-\ell} 
    = \binom{\frac{r+1}{2}-\ell-\left\lfloor\frac{\ell+1}{2}\right\rfloor}
    {2i+1-\ell}.
  \]
We now show that the indices of the last odd entries in the remaining columns of $R$ are all distinct. More precisely, we show that these indices are
\[
\ell_i \coloneq  
\begin{cases}
    2i+1 &\text{for $0\leq i<\left\lfloor\frac{r+3}{6}\right\rfloor$},\\[2pt]
    r-1-4i &\text{for $\left\lfloor\frac{r+3}{6}\right\rfloor\leq i < \left\lfloor\frac{r+3}{4}\right\rfloor$}.
\end{cases}
\]
In both cases, the binomial coefficient equals $1$, which follows in the second case from $\ell_i$ being even. Since $r$ is odd, the the indices $\ell_i$ are distinct. It remains to show that $0\leq\ell_i\leq\left\lfloor\frac{r+1}{3}\right\rfloor$ and that subsequent entries in the same column vanish.

For the first case, we observe that $2i+1\leq \left\lfloor\frac{r+1}{3}\right\rfloor$ if and only if $i<\left\lfloor\frac{r+3}{6}\right\rfloor$, and that the binomial coefficient vanishes when replacing $\ell_i$ with any larger value.

\goodbreak
For the second case, we check the extremal values of $i$. On the one hand, for
  $i = \left\lfloor\frac{r+3}{6}\right\rfloor$, we require
  \[
    \ell_i =  r-1-4 \left\lfloor\frac{r+3}{6}\right\rfloor \leq
    \left\lfloor\frac{r+1}{3}\right\rfloor.
  \]
  Indeed, this is the case for all three relevant congruence classes
  modulo $6$, namely, for $r+1 = 6s$, $r+1= 6s+2$, and $r+1=6s+4$.  On the other
  hand, for
  $i = \left\lfloor\frac{r+3}{4}\right\rfloor-1 =
  \left\lfloor\frac{r-1}{4}\right\rfloor$, we require
  \[
    \ell_i = r-1-4 \Big\lfloor\frac{r-1}{4}\Big\rfloor\geq
    0,
  \]
  which is obviously true. It remains to show that all subsequent
  entries in the same column vanish. This is certainly the case if, for $j > 0$,
  \[
    \frac{r+1}{2}-(\ell_i+j)-\left\lfloor\frac{\ell_i+j+1}{2}\right\rfloor 
    < 2i+1-(\ell_i+j).
  \]
  Using $\ell_i = r-1-4i$, we obtain the condition
  \[
    \frac{r+1}{2}-\left\lfloor\frac{r-1-4i+j+1}{2}\right\rfloor 
    < 2i+1,
  \]
  which simplifies to 
  \[
    -\left\lfloor\frac{j-1}{2}\right\rfloor < 1.
  \]
  This concludes the proof.
\end{proof}

\section{Exponential valuations}
\label{sec:expo}

This section contains the proof of the second case of \cref{thm:exp_val} and the third case of \cref{thm:tensor_val}. Let $\oZb\in\bVal$. For $P\in\PZ$, we can write 
\begin{equation}\label{eq_exdecomp}
\oZb(P) = \sum_{r\geq 0} \oZ^r(P),
\end{equation}
where $\oZ^r\colon\PZ\to\ring_r$ is a $\GL_2(\ZZ)$ equivariant valuation for each $r$. We call $\oZ^r$ the $r$-summand of $\oZb$. We start with two simple observations, which use the notation from \cref{eq_exdecomp}.

\begin{lemma}
    \label{lemma:poly_exp}
    The valuation $\oZb\in\bVal$ is translatively exponential if and only if every $\oZ^r$  is translatively polynomial with associated functions $\oZ^{r-j}$  for $0\leq j\leq r$.
\end{lemma}

\begin{proof}
    If $\oZb$ is translatively exponential, then 
    \[
    \sum_{r\geq 0}  \oZ^r(P+v) = \oZb(P+v) = \e^v \oZb(P) =\sum_{k\geq 0}\frac{v^k}{k!}\,\sum_{r\geq 0} \oZ^r(P) = \sum_{r\geq 0}\sum_{j=0}^r \oZ^{r-j}(P)\frac{v^j}{j!}
    \]
    for every $P\in\PZ$ and $v\in\ZZ^2$.
    Comparing the $r$-summands on both sides shows that $\oZ^r$ is translatively polynomial with associated valuations $\oZ^{r-j}$ for $0\leq j\leq r$. The reverse implication can be seen in the same way.
\end{proof}

\goodbreak
\begin{lemma}
    \label{lemma:hom_dilative}
    The valuation $\oZb\in\bVal$ is $d$-dilative if and only if the valuation $\oZ^r$ is $(r-d)$-homogeneous for every $r\ge0$.
    In particular, we have
    \begin{enumerate}
        \item $\bVal_d = 0$ for $d< -2$, and
        \item if $\,\oZb$ is $d$-dilative, then $\oZ^r = 0$ for every $r<d$.
    \end{enumerate}
\end{lemma}

\begin{proof}
    If $\oZb$ is $d$-dilative, then
    \begin{align*}
       \sum_{r\geq 0}\oZ^r(mP)(x,y)&= \overline\oZ(mP)(x,y) \\
       &= m^{-d}\sum_{r\geq 0}\oZ^r(P)(mx,my) = \sum_{r\geq 0}m^{r-d}\oZ^r(P)(x,y)
    \end{align*}
    for $m\in\ZZ_{>0}$.
    Comparing the $r$-summands on both sides shows that $\oZ^r$ is $(r-d)$-homogeneous. The reverse implication can be seen in the same way. The last two statements immediately follow from \cref{thm:khovanskii}.
\end{proof}
We deduce the decomposition of $\bVal$ claimed in \cref{thm:exp_val}.

\begin{proposition}
    \label{prop:exp_decomp}
    The map
    \[
    \prod_{d\geq -2}\bVal_{d} \to \bVal,\quad (\oZb_{d})_{d\geq -2} \mapsto \sum_{d\geq -2} \oZb_d
    \]
    is a linear isomorphism.
\end{proposition}

\begin{proof}
    The linearity is clear. To see that the map is injective, suppose that 
    \[\sum\nolimits_{d\geq -2}\overline\oZ_d = 0.\]
    We write $\oZb_d = \sum_{r\geq 0} \oZ^r_{r-d}$. Note that the difference of upper and lower indices is $d$, so there is no conflict of notation for different $d$. Extracting the $r$-summand of the sum, we see that
    $\sum_{d\geq -2} \oZ^r_{r-d} = 0$. It follows from \cref{lemma:hom_dilative} that $\oZ_{r-d}^r$ is $(r-d)$-homogeneous. This implies $\oZ^r_{r-d}=0$ for all $d\geq -2$. Since $r$ was arbitrary, we obtain $\overline\oZ_d=0$ for all $d\geq -2$.

    For the surjectivity, let $\overline\oZ \in\bVal$ be written as in \cref{eq_exdecomp}. By \cref{thm:khovanskii}, each $\oZ^r$ has a homogeneous decomposition $\oZ^r = \sum_{i=0}^{r+2}\oZ^r_i$ with $\oZ_i^r\in\Val_i^r$.  Since $\oZ$ is translatively exponential, \cref{lemma:poly_exp} implies that $\oZ^r$ is translatively polynomial with associated valuation $\oZ^{r-j}$ for $0\leq j\leq r$. Thus, it follows from \cref{lemma:crazywindices} that $\oZ^r_i$ is translatively polynomial with associated valuations $\oZ^{r-j}_{i-j}$ for $0\leq j\leq r$. Since the difference between rank and homogeneity degree is the same for every $j$, the reverse implication of \cref{lemma:poly_exp} shows that  the valuation 
    \[\oZb_d \coloneq  \sum\nolimits_{r\geq 0}\oZ^r_{r-d}\] 
    is translatively exponential for any $d\geq -2$. By \cref{lemma:hom_dilative}, it is  $d$-dilative. Since the $r$-summands of $\overline{\oZ}_d$ are zero for $r<d$ by \cref{lemma:hom_dilative}, in the sum $\sum_{d\geq -2} \overline{\oZ}_d(P)$, we see only finitely many polynomial terms of each degree of homogeneity. Therefore, $\sum_{d\geq -2}\overline{\oZ}_d(P)$ defines a formal power series which is equal to $\overline{\oZ}(P)$.
\end{proof}

For the results in this section, we did not use that our valuations are defined on polygons. The corresponding results can be proven similarly in the $n$-dimensional setting.

\section{Unbounded polyhedra}
\label{sec:lawrence}

Next, we classify $d$-dilative, translatively exponential valuations. Here, we are only concerned with the second case of \cref{thm:exp_val} as the remaining cases will be treated in the next section. So, for the remainder of the section, let $d>0$ be an even integer.

It follows from \cref{lemma:poly_exp} and \cref{lemma:hom_dilative} that we have 
well-defined linear maps
\begin{equation}
    \label{eq:lambda}
    \Lambda^r \colon \bVal_d \to \Val^r_{r-d},\quad \overline{\oZ} =\sum\nolimits_{s\geq 0}\oZ^s \mapsto \oZ^r.
\end{equation}
We require the following result.

\begin{lemma}
    \label{lemma:lambda_inj}
    Let $d>0$ be even. The map $\Lambda^r$ defined in \cref{eq:lambda} is injective for $r>d$ and identically zero for $r\leq d$.
\end{lemma}

\begin{proof}
    The fact that $\Lambda^r=0$ for $r\leq d$ follows from \cref{thm:khovanskii} and \cref{lemma:constterm}. To see that $\Lambda^r$ is injective for $r>d$, consider $\overline{\oZ}=\sum_{s\geq 0} \oZ^s\in\bVal_d$ with $0=\Lambda^r\oZb = \oZ^r$. It follows from \cref{lemma:poly_exp} that $\oZ^s = 0$ for all $s\leq r$, since these are the associated valuations of $\oZ^r$. Again by \cref{lemma:poly_exp}, the associated valuations of $\oZ^{r+1}\in\Val_{r-d+1}^{r+1}$ (other than $\oZ^{r+1}$ itself) are all zero. Since by assumption $r-d+1\geq 2$, \cref{lemma:familybusiness} applies and we obtain $\oZ^{r+1}=0$. Using induction on $j$, we see that $\oZ^{r+j}=0$ for all $j>0$. Thus, $\overline\oZ=0$, as desired.
\end{proof}

In the following, we want to prove that $\Lambda^{d+2}$ is a linear isomorphism.
To this end, we construct a valuation on the set of rational polyhedra $\polyhedra(\QQ^2)$ with values in the vector space 
\[\mero\coloneq  \spann\{\e^v\cdot f: v\in\RR^2,~ f\in\rats \},\]
a subspace of the space of meromorphic functions. Here, $\rats$ denotes the space of rational functions in $x$ and $y$.
The restriction of this valuation to lattice polygons will take values in $\fps$. Here, a polyhedron $Q\subseteq\RR^2$ is called \Dfn{rational} if it can be represented as $Q=\{v\in\RR^2\colon z_i(v)\leq b_i,~1\leq i \leq m\}$ for certain $z_1,\dots,z_m\in(\QQ^2)^\star$ and $b_1\dots,b_m\in\QQ$. For $X\subseteq\RR^2$, let $\pos (X)$ denote the positive hull of $X$.  We call a rational polyhedron $Q\in\polyhedra(\QQ^2)$ a \Dfn{rational cone} if $\pos(Q)=Q$, and we denote the set of rational cones by $\CZ$.

\goodbreak
We will use the following special case of \cite[Theorem 1]{lawrence}. Let $\vertices (Q)$ be the set of vertices of $Q\in\polyhedra(\QQ^2)$.

\begin{theorem}[Lawrence]
    \label{thm:lawrence}
    A valuation $\zeta_0\colon\CZ\to\mero$ that vanishes on cones containing a line can be extended to a valuation $\zeta$ on $\polyhedra(\QQ^2)$ by
    \begin{equation*}
        \label{eq:lawrence}
       \zeta(Q) \coloneq  \sum\nolimits_{v\in\vertices(Q)} \e^v\,\zeta_0(\fcone(v;Q)),
     \end{equation*}
    where $\fcone(v;Q) \coloneq  \pos(Q-v)$.
\end{theorem}

To construct our particular $\zeta_0$, we start with a valuation $\oZ_2^{d+2}\in\Val_2^{d+2}$ and consider the rational function
\begin{equation}
    \label{eq:r}
    R(x,y) \coloneq  \frac{\oZ_2^{d+2}([0,1]^2)}{xy}.
\end{equation}
Note that $R\in\rats_d$,  the space of $d$-homogeneous rational functions in $x$ and $y$. 

We require some results and notions related to cones (see, for example, \cite{BrunsGubeladze}). A pointed two-dimensional cone $C\in\CZ$ is called \Dfn{unimodular} if $C=\phi_C(\quadrant)$ with $\phi_C\in\GL_2(\ZZ)$. Here $\quadrant\coloneq \{(x,y)\in\RR^2: x,y\ge0\}$. Note that $\phi_C$ is unique up to a permutation of the coordinates. Since $R$ is symmetric in $x$ and $y$, the expression $R\circ\phi_C^\star$ is independent of the choice of $\phi_C$.
For a general two-dimensional pointed cone $K\in\CZ$, let $u_1,\dots,u_m$ be the primitive lattice vectors on the boundary of the convex hull of $K\cap\ZZ^2\setminus\{0\}$ (labeled in any direction along the boundary), where a lattice vector is \Dfn{primitive} if its coordinates are co-prime.
For a two-dimensional cone, the \Dfn{Hilbert decomposition} $\can(K)$ can be defined as 
\[\can(K) \coloneq  \{\pos\{u_{i-1},u_i\} \colon 1\leq i \leq m\}.\] 
It is easy to see that the decomposition $\can(K)$ is \Dfn{unimodular}; that is, the cones of $\can(K)$ are unimodular. 
Moreover,  
any uni\-modular decomposition of $\can(K)$ is a refinement of $\can(K)$ (cf.\ \cite[Proposition 1.19]{Oda}). 
We can describe the refinements in a structured way. For a uni\-modular cone $C=\pos\{u,v\}$, where $u,v\in\ZZ^2$ are primitive vectors, we define the \Dfn{balanced decomposition} of $C$ as $\{\pos\{u,u+v\},\pos\{u+v,v\}\}$. 
For two unimodular decompositions $\cS$ and $\cT$ of the cone $K\in\CZ$, we say that $\cT$ is an \Dfn{elementary refinement} of $\cS$ if $\cT$ can be obtained from $\cS$ by replacing a unimodular cone of $\cS$ by its balanced decomposition.
\begin{lemma}
    \label{lemma:refinements}
    Let $K\in\CZ$ be two-dimensional and pointed. For every unimodular decomposition $\cT$ of $K$, there exist unimodular decompositions $\cT_0,\dots,\cT_k$ (for some $k\in\ZZ_{\geq 0}$) of $K$ such that $\cT_0=\can(K)$, $\cT_k=\cT$, and $\cT_i$ is an elementary refinement of $\cT_{i-1}$ for all $1\leq i \leq k$.
\end{lemma}

\begin{proof}
    Since any unimodular triangulation of $K$ is a refinement of $\can(K)$, it suffices to prove the lemma for $K$ being itself a unimodular cone, that is, $\can(K)=\{K\}$ and $K=\pos\{u,v\}$, where $u$ and $v$ generate $\ZZ^2$. 
    We use induction on $m\coloneq |\cT|$. 
    For $m=1$, there is nothing to prove. So, let $m>1$ and suppose that the statement of the lemma is true for any unimodular cone $K$ and any decomposition $\cT$ of size less than $m$.    
    
    First, we show that the ray $\pos\{u+v\}$ is contained in the boundary of two cones of $\cT$. 
    To see this, we apply a unimodular transformation so that $K=\quadrant$ and, thus, $u+v = (1,1)$. 
    If $\pos\{(1,1)\}$ is not a boundary ray of a cone in $\cT$, there exists a unimodular cone $C\in \CZ$ with primitive generators $(a,b), (c,d)\in\ZZ^2$  such that $a>b\geq 0$ and $d>c\geq 0$ and such that $(1,1)$ is contained in the interior $C$. Since $C$ is unimodular, there exist $k,\ell\in\ZZ_{> 0}$ with $(1,1) = k(a,b) + \ell(c,d)$. Using that the matrix $\begin{psmallmatrix}
        a & b\\
        c & d
    \end{psmallmatrix}$ is in $\SL_2(\ZZ)$ and solving for $k$ and $\ell$ gives $k=d-c\geq 1$ and $\ell = a-b\geq 1$. This implies $(a,b)=(1,0)$ and $(c,d)=(0,1)$. Hence, $C=K$, which contradicts $m>1$. Thus, the ray $\pos\{u+v\}$ must be contained in the boundary of a cone of $\cT$. 
    
    Consequently, we obtain that $\cT$ is the disjoint union of the two unimodular decompositions 
    $\{T\in\cT\colon T\subseteq \pos\{u,u+v\}\}$ and 
    $ \{T\in\cT\colon T\subseteq \pos\{u+v,v\}\}$ of the unimodular cones $\pos\{u,u+v\}$ and $\pos\{u+v,v\}$, respectively. 
    Since each of the two decompositions contains fewer cones than $\triangulation$, the induction hypothesis applies, and we see that $\cT$ is obtained by taking the balanced decomposition of $K$ first and then refining the two halves inductively. 
\end{proof}

For a pointed two-dimensional cone $K\in\CZ$, set
\begin{equation}
\label{eq:zetanull}
\zeta_0(K) \coloneq  \sum_{C\in\can(K)} R\circ\phi_C^\star.
\end{equation}
If $K$ is not pointed or if $K$ is lower-dimensional, we set $\zeta_0(K)\coloneq 0$.

\begin{lemma}
    \label{lemma:zetanull}
    For $d>0$ even, the function $\zeta_0\colon\CZ \to \rats_d$,  defined in \cref{eq:zetanull}, is a simple, $\GL_2(\ZZ)$ equivariant valuation that vanishes on cones with lines.
\end{lemma}

\begin{proof}
    We only have to show that $\zeta_0$ is a valuation.
    Let $K\in \CZ$ and $L\subset\RR^2$ be a rational one-dimensional subspace ($L$ is generated by a rational vector) that intersects $K$. Let $K_+$ and $K_-$ denote the intersections of $K$ with the two closed halfspaces determined by $L$. By \cite[Lemma 19.5]{Barvinok2008}, it is sufficient to check that
    \begin{equation}
    \label{eq:weak_val}
    \zeta_0(K) = \zeta_0(K_+) + \zeta_0(K_-)
    \end{equation}
    for all $K\in\CZ$ and rational one-dimensional subspaces $L$.
    
    If $\dim K<2$, this is trivial. So, let $K$ be two-dimensional. 
    First, we assume that $K$ is pointed. Let $\can\coloneq \can(K)$ and $\can_\sigma \coloneq  \can(K_\sigma)$ for $\sigma \in \{-,+\}$. Clearly, $\cT \coloneq  \can_- \cup \can_+$ is a unimodular decomposition of $K$. As such, it is obtained from $\can$ by a finite sequence of elementary refinements using\cref{lemma:refinements}. Consider a unimodular cone $C\in\cT$. By $\GL_2(\ZZ)$ equivariance, we may assume that $C=\quadrant$. If we  show that 
    \[\zeta_0(\quadrant)=\zeta_0(\pos\{e_1,e_1+e_2\}) + \zeta_0(\pos\{e_2+e_1,e_2\}),\]
    then we can use the $\GL_2(\ZZ)$ equivariance and apply \cref{lemma:refinements} inductively to obtain \cref{eq:weak_val}. Hence, by definition of $\zeta_0$, we need to verify that
    \begin{equation}
    \label{eq:func_eq_1}
    R(x,y) = R(x,x+y) + R(x+y,y)
    \end{equation}
    for $x,y\in\RR$. Writing $P_0 \coloneq  [0,1]^2$, $P_1\coloneq \begin{psmallmatrix}
         1 & 1\\
         0 & 1       
    \end{psmallmatrix}P_0$, and $P_2 \coloneq  \begin{psmallmatrix}
         1 & 0\\
         1 & 1       
    \end{psmallmatrix}P_0$, we have
    \[
    \begin{split}
        R(x,x+y) = \frac{\oZ_2^{d+2}(P_0)(x,x+y)}{x(x+y)} = \frac{\oZ_2^{d+2}(P_1)}{x(x+y)}.
    \end{split}
    \]
    Similarly, we obtain
    \[R(x+y,y) = \frac{\oZ_2^{d+2}(P_2)}{(x+y)y}.\] 
    Thus,
    \[
        R(x,x+y) + R(x+y,y) = \frac{x \oZ_2^{d+2}(P_2) + y\oZ_2^{d+2}(P_1)}{xy(x+y)}.
    \]
    \goodbreak\noindent
    Let $T\coloneq \begin{psmallmatrix}
        1 & 1\\
        0 & 1
    \end{psmallmatrix}\sttriangle$. Using the decompositions \[P_0=T\sqcup (-T+e_1+e_2), \,P_1 = T\sqcup (-T+2e_1+e_2), \, P_2=(-T+e_1+e_2)\sqcup(e_2+T),\] we compute
    \begin{align*}
        x(\oZ_2^{d+2}(P_2) - \oZ_2^{d+2}(P_0)) &= x(\oZ_2^{d+2}(T+e_2) - \oZ_2^{d+2}(T)) \\
        &= xy\oZ_1^{d+1}(T)\\
        &=y(-x\oZ_1^{d+1}(-T+e_1+e_2))  \\
        &= y(\oZ_2^{d+2}(-T+e_1+e_2) - \oZ_2^{d+2}(-T+2e_1+e_2))\\
        &=y(\oZ_2^{d+2}(P_0) - \oZ_2^{d+2}(P_1)).
    \end{align*}
    We have used that $\oZ_2^{d+2}$ is simple by \cref{lemma:simple} and that $\oZ_1^{d+1}$ is odd and translation invariant by \cref{lemma:onehom_translation_invariant}.
    Hence, 
    \[x\oZ_2^{d+2}(P_2) + y\oZ_2^{d+2}(P_1) = (x+y)\oZ_2^{d+2}(P_0)\]
    and \cref{eq:func_eq_1} follows. So we obtain \cref{eq:weak_val} if $K$ is a pointed unimodular cone. 

    \goodbreak
    Next, let $K$ be a halfspace. We may assume that $L$ intersects the interior of $K$. Then, $K_+$ is unimodular if and only if $K_-$ is unimodular. First, assume that the cones $K_+$ and $K_-$ are unimodular. By $\GL_2(\ZZ)$ equivariance, it suffices to consider the case where $K$ is the upper halfplane and $L$ is the $y$-axis. In this case, \cref{eq:weak_val} is equivalent to the functional equation
    \begin{equation}
        \label{eq:func_eq2}
        0 = R(x,y) + R(-x,y)
    \end{equation}
    for $x,y\in\RR$. Since $\oZ_2^{d+2}([0,1]^2)$ is  $\sqgroup$ invariant, \cref{eq:func_eq2} is  fulfilled.
    
    Next, we assume that $K_+$ is a non-unimodular cone. We can choose a one-dimensional rational subspace $L^\prime$ that meets the interior of $K_+$ such that the cone $K_+^\prime$ bounded by the boundary of $K$ and $L^\prime$ is unimodular. Let $K_{-}^\prime$ be the closure of $K\setminus K_+^\prime$, which is a unimodular cone. Thus,
    \[
        \zeta_0(K) = \zeta_0(K_+^\prime) + \zeta_0(K_-^\prime).
    \]
    Moreover, $K_{-}^\prime$ is decomposed by the subspace $L$ into the cones $K_-$ and $K_+^{\prime\prime}$, where $K_+^{\prime\prime}$ is the closure of $K_+\setminus K_+^\prime$. Since we already proved that \cref{eq:weak_val} holds for the pointed cone $K_-^\prime$ we deduce
    \[
        \zeta_0(K) = \zeta_0(K_+^\prime) + \zeta_0(K_+^{\prime\prime}) + \zeta_0(K_-) = \zeta_0(K_+) + \zeta_0(K_-), 
    \]
    which verifies \cref{eq:weak_val} for halfspaces. Since the statement is trivial for $K=\RR^2$, the proof is complete.
\end{proof}

Using \cref{thm:lawrence}, we obtain a simple, $\GL_2(\ZZ)$ equivariant, and translatively exponential  valuation
\begin{equation}
    \label{eq:zeta}
    \zeta\colon\polyhedra(\QQ^2) \to \mero, \quad \zeta(P) \coloneq  \sum\nolimits_{v\in\vertices(P)} \e^v\zeta_0(\fcone(v;P)).
\end{equation}
Moreover, it follows from $\zeta_0(\fcone(v;P))\in\rats_d$ that $\zeta$ is $d$-dilative. Let us evaluate $\zeta$ at the standard triangle $\sttriangle$.

\goodbreak
\begin{lemma}
    \label{lemma:zeta_st}
    For $\oZ_2^{d+2}\in\Val_2^{d+2}$ with $d>0$ even, the valuation $\zeta$ defined in \cref{eq:zeta} satisfies
    \[\begin{split}
        \zeta(\sttriangle) &= 2\oZ_2^{d+2}(\sttriangle)\cdot\left(\frac{\e^y-1}{y(y-x)} - \frac{\e^x-1}{x(y-x)}\right) - \oZ_1^{d+1}(\sttriangle)\cdot\left( \frac{x(\e^y-1)}{y(y-x)} - \frac{y(\e^x-1)}{x(y-x)} \right)\\
        &= 2\oZ_2^{d+2}(\sttriangle) \sum_{m > 0}\frac{[m]_{x,y}}{(m+1)!} - \oZ_1^{d+1}(\sttriangle)\sum_{m>0} \frac{[m]_{x,y}-x^{m-1}-y^{m-1}}{m!},
    \end{split}
    \]
    where $[m]_{x,y} \coloneq  \sum_{j=0}^{m-1}x^jy^{m-j}$ and $\oZ_1^{d+1}\in\Val_1^{d+1}$ is the associated valuation of $\oZ_2^{d+2}$. In particular, the function $\zeta(\sttriangle)$ is analytic. 
\end{lemma}

\begin{proof}
    We set $f_2: = \oZ_2^{d+2}(\sttriangle)$ and $f_1\coloneq \oZ_1^{d+1}(\sttriangle)$. Since $\oZ_{i}^{d+i}$ is simple for $i\in\{1,2\}$ by \cref{lemma:simple}, considering the triangulation $[0,1]^2 = \sttriangle \sqcup (-\sttriangle+e_1+e_2)$, we obtain
    \begin{equation}
        \label{eq:r_intermsof_f}
        R(x,y) = \frac{2f_2(x,y) - (x+y)f_1(x,y)}{xy}.
    \end{equation}
    For the standard triangle, we have \[
    \fcone(0;\sttriangle) = \quadrant,~\fcone(e_1;\sttriangle) = \phi_1\quadrant,~\text{and}~\fcone(e_2;\sttriangle) = \phi_2\quadrant,
    \]
    where $\phi_1 = \begin{psmallmatrix}
        -1 & -1\\
        \phantom{-}1 & \phantom{-}0
    \end{psmallmatrix}$ and
    $\phi_2 = \begin{psmallmatrix}
        \phantom{-}0 & \phantom{-}1\\
        -1 & -1
    \end{psmallmatrix}$. Note that $\phi_i\sttriangle = \sttriangle-e_i$ for $i\in\{1,2\}$. So, in particular, we have $\phi_1,\phi_2\in\group$.

    Using the representation \cref{eq:r_intermsof_f} of $R$ and the $\group$ invariance of $f_1$, we compute
    \[
    \begin{split}
        R\circ \phi_2^\star &= \frac{1}{(x-y)(-y)}\left( 2f_2\circ\phi_2^\star - (x-y + (-y))f_1\circ\phi_2^\star \right)\\
        &=\frac{1}{y(y-x)}\left( 2\oZ_2^{d+2}(\sttriangle-e_2) - (x-2y)f_1  \right)\\
        &=\frac{1}{y(y-x)}\left( 2(f_2 - yf_1) -(x-2y)f_1\right)\\
        &=\frac{1}{y(y-x)}\left( 2f_2  -xf_1\right).
    \end{split}
    \]
    Following the same lines, we compute
    \[
         R\circ \phi_1^\star = \frac{1}{x(x-y)}\left( 2f_2  -yf_1\right).
    \]
    So we obtain
    \[\begin{split}
    \zeta(\sttriangle) &= \e^0\zeta_0(\RR^2_{\geq 0}) + \e^x\zeta_0(\phi_1 \RR^2_{\geq 0}) + \e^y\zeta_0(\phi_2\RR^2_{\geq0}) \\[4pt]
    &= R+\e^x\, R\circ\phi_1^\star + \e^y\,R\circ\phi_2^\star\\[4pt]
    &=\frac{2f_2-(x+y)f_1}{xy} + \e^x\,\frac{2f_2  -yf_1}{x(x-y)} + \e^y\,\frac{2f_2  -xf_1}{y(y-x)} \\
    &=2f_2\cdot\left(\frac{\e^y-1}{y(y-x)} - \frac{\e^x-1}{x(y-x)}\right) - f_1\cdot\left( \frac{x(\e^y-1)}{y(y-x)} - \frac{y(\e^x-1)}{x(y-x)} \right).
    \end{split}\]
    This proves the first identity. 
    
    For the second identity, note that $\tfrac{e^x-1}{x} = \sum_{m\geq 0} \tfrac{x^m}{(m+1)!}$ is an analytic function and that $y^m-x^m = (y-x)[m]_{x,y}$. For the first bracket term, we thus obtain
    \[
    \frac{\e^y-1}{y(y-x)} - \frac{\e^x-1}{x(y-x)} = \sum_{m\geq 0} \frac{y^m-x^m}{(m+1)!(y-x)}= \sum_{m>0} \frac{[m]_{x,y}}{(m+1)!}.
    \]
    In particular, the left-hand side is an analytic function on all $\RR^2$.
    Similarly, we have for the second bracket term:
    \[\begin{split}
     \frac{x(\e^y-1)}{y(y-x)} - \frac{y(\e^x-1)}{x(y-x)} &= \sum_{m\geq 0} \frac{xy^m-yx^m}{(m+1)!(y-x)}\\
     &= -1+\sum_{m> 0} \frac{xy(y^{m-1}-x^{m-1})}{(m+1)!(y-x)}   = -1 + \sum_{m>0}\frac{xy[m-1]_{x,y}}{(m+1)!}\\
     &= \sum_{m\geq 0}\frac{[m+1]_{x,y}-x^m-y^m}{(m+1)!} = \sum_{m>0} \frac{[m]_{x,y}-x^{m-1}-y^{m-1}}{m!}.
    \end{split}\]
    Again, the left-hand side is an analytic function on all $\RR^2$.
\end{proof}

\goodbreak
Let $\oZb$ be the restriction of $\zeta$ to $\PZ$. By \cref{lemma:zeta_st}, we have $\oZb(\sttriangle)\in\fps$.
Since $\zeta$ is simple, it follows that $\oZb\in\bVal_d$. 
From the series expansion in \cref{lemma:zeta_st}, we see that
the $(d+1)$-summand of $ \oZb$ is $\oZ_1^{d+1}$ and the $(d+2)$-summand of $\oZb$ is $\oZ_2^{d+2}$. This shows that the map $\Lambda^r$ is surjective. Combined with \cref{lemma:lambda_inj}, we obtain the following result.
\begin{lemma}
    \label{lemma:lambda_iso}
    For $d>0$ even, the map $\Lambda^{d+2} \colon \bVal_{d+2} \to \Val^{d+2}_{2}$  defined in \cref{eq:lambda} is a linear isomorphism.
\end{lemma}

\goodbreak
This allows us to compute the dimensions of several spaces of equivariant valuations.

\begin{corollary}
    \label{cor:d_even_dim}
    Let $d>0$ be even. Then we have
    \[
        \dim\bVal_d = \dim \Val_i^{d+i} = \parts{\tfrac d2 + 1}
    \]
    for all $i\geq 2$.
\end{corollary}

\begin{proof}
    By \cref{thm:2hom}, we have $\dim\Val_2^{d+2} = \parts{\tfrac d2 + 1}$. By \cref{lemma:lambda_iso}, we have $\dim\bVal_d = \parts{\tfrac d2 + 1}$. Now let $i\geq 2$ and let 
    $A_i^{d+i} \colon \Val_i^{d+i} \to \Val_{i-1}^{d+i-1}$ be the linear operator that maps a valuation $\oZ_i^{d+i}$ to its associated valuation of rank $d+i-1$. 
    Note that this associated valuation is $(i-1)$-homogeneous by \cref{lemma:crazywindices}. Since $i\geq 2$, applying \cref{lemma:familybusiness} inductively shows that each $A_i^{d+i}$ with $i\geq 2$ is injective. So, the linear map $A_i^{d+i} \circ \cdots \circ A_3^{d+3}$ from $\Val_i^{d+i}$ to $\Val_2^{d+2}$ is injective. At the same time, the map $\Lambda^{d+i}\colon \bVal_d\to\Val_i^{d+i}$ is injective by \cref{lemma:lambda_inj}. This shows the claim for all $i\geq 2$ (cf.\ \cref{fig:diagram}).
\end{proof}

\goodbreak
Let us have a closer look at the valuations we obtain. For a lattice polygon $P\in\PZ$ and a vertex $v$ of $P$, let $u_1(v),\dots,u_{m_v}(v)$ be the lattice points on the boundary of the convex hull of $P\cap\ZZ^2\setminus\{v\}$ that are visible from $v$ (labeled in any direction along the boundary). Let $\phi_i[v]$ denote the linear transformation with $\phi_i[v](e_j) = u_{i+j-2}(v)-v$ for $1\leq i\leq m_v$ and $1\leq j\leq 2$. Then we have $\phi_i[v]\in\GL_2(\ZZ)$ and it follows from \cref{eq:zeta} that
\begin{equation}
\label{eq:triang_free}
\overline\oZ(P) = \sum_{v\in\vertices(P)} \e^v\,\sum_{i=1}^{m_v} R\circ \phi_i[v]^\star.
\end{equation}
By \eqref{eq:r}, the function $R$ is determined by $\oZ_2^{d+2}$. As we saw in \cref{prop:2hom_equiv}, there is a one-to-one correspondence between valuations $\oZ_2^{d+2}\in\Val_2^{d+2}$ and $\group$~invariant polynomials $f\in\ring_{d+1}^\group$ for which there exists an $h\in\ring_{d+2}^\group$ such that $h+(x+y)f$ is $\sqgroup$ invariant. If it exists, this $h$ is uniquely determined by $f$ by \cref{prop:2hom_equiv}. Hence,  \cref{eq:triang_free} allows us to compute $\overline{\oZ}(P)$ solely in terms of the initial polynomial $f$ and the lattice geometry of $P$. In particular, no triangulation is needed.

Similarly, we obtain a triangulation-free description for the valuations in $\Val_i^{d+i}$ for $i>1$ as the $(d+i)$-summand of the series in \cref{eq:triang_free}. For $i=1$, we have this description only for valuations of the form $\oZ_f$ (as defined in \cref{eq:adhoc_definition}), where the polynomial $f\in\ring_{d+1}^\group$ is such that $(x+y)f + \ring_{r+1}^\group$ intersects $\ring_{r+1}^\sqgroup$.

\section{The remaining cases}
\label{sec:remaining}

In this section, we consider the cases in \cref{thm:tensor_val} and \cref{thm:exp_val} where the respective vector spaces are one-dimensional. The arguments used here essentially appear in \cite{ludwigsilverstein} already; we recollect them to keep the paper self-contained.

\begin{proposition}
    \label{prop:unique_ehrhart}
       Let $r\in\ZZ_{\geq 0}$ and $i\in\{1,\dots,r+2\}$. If $r\leq i\leq r+2$ or $r-i$ is odd, then $\V_i^r = \spann\{\Lat_i^r\}$. Moreover, $\Lat_i^r$ is not identically zero in these cases and $\Lat_i^r([0,e_1])\ne0$ unless $i=r+2$.
\end{proposition}

\begin{proof}
    We start by showing that $\Lat_i^r$ for $r-i$ odd is not identically zero on $\cP(\ZZ^1)$ and thus not identically zero on $\PZ$.
    Using Faulhaber's formula (see \cite{ludwigsilverstein} for more information), we obtain for $m\ge1$ that
    \[
    \Lat^r(m\,[0,e_1])(x,y) = \frac{1}{r!}\sum_{k=0}^m k^rx^r = \frac{x^r}{r!(r+1)}\sum_{\ell=0}^r (-1)^\ell \binom{r+1}{\ell}B_\ell\, m^{r+1-\ell},
    \]
    where $B_\ell$ denotes the $\ell$th Bernoulli number. Since $B_\ell\neq 0$ for $\ell$ even, this shows that $\Lat_i^r([0,e_1])\neq 0$ for $r-i$ odd.
    Since also $B_1\neq 0$, we obtain the claim that $\Lat_i^r([0,e_1])\ne0$ holds for all $i$ and $r$ as in the proposition with $i\neq r+2$. Moreover, $\Lat_{r+2}^r$ is not identically zero by \cite[Lemma 26]{ludwigsilverstein}.

    Next, for $r=0$, it is easy to see that  $\Lat^0_i$ is not identically zero on $\PZ$ for $i\in\{0,1,2\}$.
    Thus, for $r\leq i\leq r+2$, using \cref{lemma:crazywindices}, we see that $\Lat_i^r$ is not translation invariant since $\Lat_{i-r}^0$ is among its associated functions. So, in particular, $\Lat_i^r\ne0$. So, we have established $\Lat_i^r\ne0$  for the indices occurring in the proposition.

    We note that $\V_i^0 = \spann\{\Lat_i^0\}$ by the Betke--Kneser theorem. For $r\geq 2$ even, consider the one-homogeneous valuation $\oZ_1^r$, which is translation invariant by \cref{lemma:onehom_translation_invariant}. Define the valuation $\oZ^\prime\colon\cP(\ZZ^1)\to\ring_r$ by $\oZ^\prime(S) \coloneq  \oZ_1^r(S\times\{0\})$. By \cite[Lemma 5]{ludwigsilverstein}, we have $\oZ^\prime(S)\in\RR[x]_r$. It is elementary to check that there is only one even, one-homogeneous, translatively polynomial valuation from $\cP(\ZZ^1)$ to $\RR[x]_r$. Thus, $\oZ^\prime = c\Lat_1^r$ for some $c\in\RR$. It follows that the valuation $\oY\coloneq \oZ_1^r-c\Lat_1^r$ is simple. Lemma 24 from \cite{ludwigsilverstein} yields that $\oY([0,1]^2)=0$. Since we have $[0,1]^2 = \sttriangle \sqcup(-\sttriangle+e_1+e_2)$, it follows from the assumption that $r$ is even, and the fact that $\oY$ is simple, translation invariant, and $\GL_2(\ZZ)$ equivariant that $\oY(\sttriangle)=0$. But this also implies $\oY=0$. So, we have established the statement for one-homogeneous valuations of even rank.

    Next, suppose that $i=r=1$. Then, $\Lat_1^1$ is called the \emph{discrete Steiner point} and it follows from Theorem 5 of \cite{BoeroeczkyLudwig} that $\Val_1^1 = \spann\{\Lat_1^1\}$. Hence, we proved the statement for $i=1$.

    Finally, let $r>0$ and $i>1$ satisfy the proposition's assumptions. Let $\oZ_i^r\in\V_i^r$. Then, $r-1$ and $i-1$ also satisfy the assumptions. By induction on $i$, we have $\oZ_{i-1}^{r-1}=c\Lat_{i-1}^{r-1}$. Since the associated functions are hereditary by \cref{lemma:hereditary}, this means that 
    $$\oZ_{i-j}^{r-j} = c\Lat_{i-j}^{r-j}$$ 
    for all $0<j\leq\min\{r,i\}$. 
    Consequently, $\oZ_i^r-c\Lat_i^r$ is translation invariant. Using the homogeneous decomposition theorem for translation invariant valuations \cite{McMullen77}, we obtain that the polynomial $m\mapsto (\oZ_i^r -c_i\Lat_i^r)(mP)$ has degree 2. This yields $\oZ_i^r-c\Lat_i^r=0$ as desired, except for $i=2$. But then the claim follows from \cite[Proposition 23]{ludwigsilverstein}.
\end{proof}

\goodbreak
For $d \geq -2$, we set $\bLat_d \coloneq  \sum_{r\geq 0} \Lat_{r-d}^r$. By \cref{lemma:poly_exp} and \cref{lemma:hom_dilative}, we have $\bLat_d\in\bVal_d$. \cref{prop:unique_ehrhart} allows us to determine the space $\Val_d$ for $d$ odd or non-positive.

\begin{corollary}
    \label{cor:d_odd_dim}
    If $d>0$ is odd or $d\in\{-2,-1,0\}$, then $\bVal_d = \spann\{\bLat_d\}$. Moreover, $\bLat_d([0,e_1])$ is not identically zero in these cases  and $\bLat_d([0,e_1])\ne0$  unless $d=-2$.
\end{corollary}

\begin{proof}
    The fact that $\bLat_d([0,e_1])\neq 0$ for $d=0$ or $d$ odd follows from $\Lat_1^{d+1}([0,e_1])\neq 0$.
    To see that $\bLat_{-2}$ is not identically zero, it suffices to note that its constant part is $\Lat_2^0$, the two-dimensional volume.
    
    Let $\oZb\in\bVal_d$ and consider for $r>d$ its $r$-summand $\oZ^r$. By \cref{lemma:hom_dilative}, we have $\oZ^r \in \Val_{r-d}^r$. By the assumption on $d$, \cref{prop:unique_ehrhart} applies, so there exists a constant $c_r\in\RR$ with $\oZ^r = c_r\Lat_{r-d}^r$. By \cref{lemma:poly_exp}, the valuations $\oZ^r$ for $r\geq 0$ are associated functions of one another. Since also $\Lat_{r-d}^r$ for $r\geq 0$ are associated functions of one another, we see that $c_r = c$ for all $r\geq 0$ and some $c\in\RR$ that does not depend on $r$. Thus, $\oZb = c \bLat_d$, as claimed.
\end{proof}

\section{Proofs of the main results}
\label{sec:proofs}

Here, we summarize how the results of the previous sections imply the main results \cref{thm:tensor_val} and \cref{thm:exp_val}.

\begin{proof}[Proof of \cref{thm:tensor_val}]
The first case follows from \cref{prop:unique_ehrhart} for $i\geq 1$, and from \cref{thm:betke-kneser} for $i=r=0$. The second case is \cref{cor:onehom_dimension}. The third case is \cref{cor:d_even_dim}. The fourth case consists of the two subcases, $i=0<r$ and $i>r+2$. In the former case, the statement follows from \cref{lemma:constterm}, and in the latter case, it follows from \cref{thm:khovanskii}.
\end{proof}

\begin{proof}[Proof of \cref{thm:basis}]
    The polynomials $\ff_{k,\ell} \coloneq  \ptwo^k\pthree^\ell$ with $2k+3\ell=r$ form a basis of $\ring_r^\group$ by \cref{lemma:chevalley}. By \cref{thm:eval_iso}, the evaluation map $\Theta^r$ is an isomorphism from $\Val_1^r$ to $\ring^\group_r$, if $r>1$ is odd. Hence, 
    \[
    \begin{split}
        \Val_1^r &= \spann\{ (\Theta^r)^{-1}(\ff_{k,\ell}) \colon 2k+3\ell=r\}\\
                 &= \spann\{ \oZ_{\ff_{k,\ell}} \colon 2k+3\ell=r\},   
    \end{split}
    \]
    where $\oZ_{\ff_{k,\ell}}$ is given by \cref{eq:adhoc_definition}. Note that 
    that
    $
    \Lat^r_1(\sttriangle)  \in \ring^\group_r
    $
    also if $r$ is even.
    From \cref{lemma:chevalley}, we obtain $\ring_r^\group = \spann\{\mathbf{p}_r\}$ for $r\in\{2,3\}$.
    Hence, $\ff_{k,\ell}$ is a multiple of $\Lat_1^2(\sttriangle)^k\Lat_1^3(\sttriangle)^\ell$. By \cite[Lemma 28]{ludwigsilverstein}, this multiple is non-zero, which concludes the proof.
\end{proof}

\begin{proof}[Proof of \cref{thm:exp_val}]
    The statement on the dilative decomposition is proven in \cref{prop:exp_decomp}.
    For the determination of dimensions, the first case is proved in \cref{cor:d_odd_dim}, the second case in \cref{cor:d_even_dim}, and the third case in \cref{lemma:hom_dilative}.
\end{proof}

\goodbreak
\section{Consequences}
\label{sec:consequences}

We collect some consequences of our results.
\subsection*{Minkowski additive functions}
\cref{thm:basis} can be rewritten as a classification of Minkowski additive functions on $\PZ$. A function  $\oZ$ on $\PZ$ with values in a vector space is \Dfn{Minkowski additive} if
\[\oZ(P+Q)=\oZ(P)+\oZ(Q)\] 
for all $P,Q\in\PZ$, where $P+Q=\{v+w\colon v\in P,\, w\in Q\}$ is the Minkowski sum of $P$ and $Q$. Since \[P+Q=(P\cup Q)+(P\cap Q)\] 
for $P,Q\in\PZ$ such that also $P\cap Q, P\cup Q\in\PZ$, every Minkowski additive function is a valuation, which is easily seen to be one-homogeneous. We remark that every translation invariant, one-homogeneous valuation is Minkowski additive (cf.\ \cite[Theorem 6]{McMullen77}).

\begin{corollary}
A function $\oZ\colon\PZ\to \ring$ is translation invariant, $\glz$ equivariant, and Minkowski additive if and only if there are $m\in\ZZ_{\ge 0}$ and constants $c_r, c_{k,\ell}\in\RR$  with $r,k,\ell\in\ZZ_{\ge0}$ such that
\begin{equation*}
\oZ(P)= \sum_{r=0}^m c_{2r} \Lat_1^{2r}(P)+ \sum_{r=0}^3 c_{2r+1} \Lat_1^{2r+1}(P)
    +\sum_{r=4}^m \sum_{\,\,2k+3\ell=2r+1} c_{k,\ell} \Lat_1^{2k,3\ell}(P)  
\end{equation*}
for every $P\in\PZ$. 
\end{corollary}

\noindent
Here, we have used again that $\parts{r}=1$ for $r\in\{3,5,7\}$.

\subsection*{Simple, translatively exponential valuations}

By \cref{lemma:simple}, each valuation in $\Val_i^{d+i}$ is simple, if $d>0$ is even. Hence, \cref{lemma:poly_exp} implies that each valuation in $\bVal_d$ is also simple. In fact, for each $i$, the $i$-summand $\oZ_i$ of a valuation $\oZb\in\bVal_d$ is in $\Val_i^{d+i}$ and thus simple by \cref{eq:lambda}.

If $d\geq -2$ is \emph{not} a positive and even number, then, by \cref{cor:d_odd_dim}, we have $\bVal_d = \spann\{\bLat_d\}$. If $d=-2$, the valuation $\bLat_{(-2)}$ is, by definition, the sum of the top-degree components $\Lat_{r+2}^r$ in \cref{eq:EC}. Since these correspond to the (continuous) moments by \cite[Lemma 26]{ludwigsilverstein}, we have
\begin{equation}
\label{eq:laplace}
\bLat_{(-2)}(P) = \sum_{r\geq 0} \Lat_{r+2}^r(P) = \sum_{r\geq 0} \frac{1}{r!}\int_P v^r \mathrm{d} v = \int_P \e^v\mathrm d v.
\end{equation}
This valuation is called \emph{the} exponential valuation in \cite{Barvinok2008}. We introduced it as $\expval(P)$ in the introduction. It is the Laplace transform of the indicator function of $-P$ (see \cite{LiMa}). For our purposes, it is interesting to observe that \cref{eq:laplace} implies that $\bLat_{(-2)}$ is simple. 

For $d>-2$ odd and for $d=0$, the valuations $\bLat_d$ are not simple by \cref{cor:d_odd_dim}. Since also their first non-zero $r$-summand is not simple by \cref{prop:unique_ehrhart}, it follows that no sum of the valuations $\bLat_d$, where $d$ is odd or zero, is simple. Hence, we obtain the following characterization of simple, translatively exponential, $\glz$ equivariant valuations from \cref{thm:exp_val}.

\begin{corollary}
    \label{cor:simple_expval}
    We have
    \[\{\oZb\in \bVal: \oZb \text{ is simple}\} = \prod_{\mathclap{\substack{d\ge -2,\\ d \text{ even},\, d\ne 0}}} \,\,\,\bVal_{d}.\]
\end{corollary} 

\noindent
Thus, the space of simple, translatively exponential, $\glz$ equivariant valuations consists of the valuations we constructed in \cref{sec:lawrence} together with $\bLat_{(-2)}$. However, one should not regard $\bLat_{(-2)}$ as an ``exceptional'' valuation. In fact, it can be obtained in the same way as the valuations in $\bVal_d$ with $d>0$ even by choosing $R(x,y)= \tfrac{1}{xy}$ in \cref{eq:r} and using the previous construction from there on. For this particular case, it is also described in \cite[Chapter 8]{Barvinok2008} how to evaluate $\bLat_{(-2)}(P)$ based on an extension to rational polyhedra.

\smallskip
\subsection*{Acknowledgments}
A.~Freyer and M.~Ludwig were supported, in part, by the Austrian Science Fund (FWF):  10.55776/P34446. We used \texttt{SageMath} \cite{sagemath} to produce examples that motivated our results.

\bibliography{bibliography}
\bibliographystyle{acm}
\end{document}

%% file: macros.tex
\newcommand{\RR}{\mathbb{R}}
\newcommand{\ZZ}{\mathbb{Z}}

\newcommand{\QQ}{\mathbb{Q}}

\DeclareMathOperator{\spann}{span}
\DeclareMathOperator{\conv}{conv}
\DeclareMathOperator{\vol}{vol}
\DeclareMathOperator{\Lat}{L}
\DeclareMathOperator{\bLat}{\overline{L}}
\DeclareMathOperator{\GL}{GL}
\DeclareMathOperator{\SL}{SL}
\newcommand{\ring}{\RR[x,y]}
\newcommand{\fps}{\RR[[x,y]]}
\newcommand{\rats}{\RR(x,y)}

\DeclareMathOperator{\centroid}{cen}
\newcommand{\e}{\mathrm{e}}

\DeclareMathOperator{\pos}{pos}
\newcommand{\quadrant}{\RR^2_{\geq 0}}

\newcommand{\group}{\mathbf{G}}

\newcommand{\sqgroup}{\mathbf{D}}
\newcommand{\thirdgroup}{\mathbf{H}}

\renewcommand{\gg}{\mathbf{g}}
\newcommand{\ff}{\mathbf{f}}
\newcommand{\pp}{\mathbf{p}}
\newcommand{\qq}{\mathbf{q}}
\newcommand{\ptwo}{\pp_2}
\newcommand{\pthree}{\pp_3}
\def\gprd{\qq}
\def\gsqr{\tilde\qq}

\newcommand{\PZ}{\mathcal{P}(\ZZ^2)} 
\newcommand{\CZ}{\mathcal{C}(\QQ^2)} 
\newcommand{\US}{\mathcal{S}(\ZZ^2)} 

\newcommand{\polyhedra}{\mathcal{Q}}
\newcommand{\V}{\mathrm{Val}}

\DeclareMathOperator{\can}{\mathcal{H}}

\DeclareMathOperator{\fcone}{fcone}

\newcommand{\eval}{\Theta}
\newcommand{\parts}[1]{p_{2,3}\left(#1\right)}
\newcommand{\expval}{\operatorname{\overline E}}

\newcommand{\mero}{\mathcal{V}}
\newcommand{\vertices}{\operatorname{vert}}

\newcommand{\cS}{\mathcal{S}}
\newcommand{\cT}{\mathcal{T}}

\newcommand{\cP}{\mathcal{P}}

\newcommand{\Val}{\mathrm{Val}}
\newcommand{\bVal}{\overline{\Val}}

\newcommand{\sttriangle}{\mathbf T} 
\newcommand{\triangulation}{\mathcal{T}}

\newcommand{\Dfn}[1]{\emph{#1}} 

\newcommand{\oZ}{\operatorname{Z}} 
\newcommand{\oz}{\operatorname{Z}'} 
\newcommand{\oY}{\operatorname{Y}}

\newcommand{\oX}{\operatorname{X}}
\newcommand{\oZb}{\operatorname{\overline{Z}}}
\newcommand{\oZt}{\operatorname{\tilde{Z}}}

\def\lpn{{\mathcal P}(\ZZ^n)} 

\newcommand{\slnz}{\operatorname{SL}_n(\ZZ)} 
\newcommand{\glnz}{\operatorname{GL}_n(\ZZ)} 
\newcommand{\glz}{\operatorname{GL}_2(\ZZ)} 